\documentclass{amsart}

\usepackage{amssymb}
\usepackage{amscd}

\title[Geometric interpretation of double shuffle relation]
{Geometric interpretation of double shuffle relation for multiple $L$-values}

\author{Hidekazu Furusho}

\address{Graduate School of Mathematics, Nagoya University, 
Furo-cho, Chikusa-ku, Nagoya, 464-8602, Japan}
\email{furusho@math.nagoya-u.ac.jp}

\newcommand{\Exp}{\text{Exp}}
\newtheorem{thm}{Theorem}[section]
\newtheorem{lem}[thm]{Lemma}

\newtheorem{prop}[thm]{Proposition}  

\theoremstyle{remark}

\theoremstyle{definition}
\newtheorem{defn}[thm]{Definition}
\newtheorem{rem}[thm]{Remark}

\newtheorem{pf}{Proof}

\numberwithin{equation}{section}
\begin{document}
\bibliographystyle{amsalpha+}
\maketitle
\begin{abstract}
This paper gives a geometric interpretation of 
the generalized (including the regularization relation) double shuffle relation for multiple $L$-values.
Precisely it is proved that Enriquez' mixed pentagon equation implies the relations.
As a corollary, an embedding
from his cyclotomic analogue of the Grothendieck-Teichm\"{u}ller group 
into Racinet's cyclotomic double shuffle group is obtained.
It cyclotomically extends the result of our previous paper \cite{F08}
and the project of Deligne and Terasoma
which are the special case $N=1$ of our result.
\end{abstract}

\tableofcontents

\setcounter{section}{-1}
\section{Introduction}
Multiple $L$-values $L(k_1,\cdots,k_m;\zeta_1,\cdots,\zeta_m)$
are the complex numbers defined by the following series
\begin{equation}\label{multiple L-value}
L(k_1,\cdots,k_m;\zeta_1,\cdots,\zeta_m)
:=\sum_{0<n_1<\cdots<n_m}\frac{\zeta_1^{n_1}\cdots \zeta_m^{n_m}}
{n_1^{k_1}\cdots n_m^{k_m}}
\end{equation}
for $m$, $k_1$,\dots, $k_m\in {\bf N} (={\bf Z}_{>0})$
and $\zeta_1$,\dots,$\zeta_m\in\mu_N$(:the group of 
$N$-th roots of unity in $\bold C$).
They converge if and only if $(k_m,\zeta_m)\neq (1,1)$.
Multiple zeta values are regarded as a special case for
$N=1$.
These values have been discussed 
in several papers \cite{AK, BK, G, R} etc.
Multiple $L$-values appear as coefficients of
the cyclotomic Drinfel'd associator
$\varPhi^N_{KZ}$
\eqref{cyclotomic explicit formula}
in $U\frak F_{N+1}$:
the non-commutative formal power series ring with $N+1$ variables
$A$ and $B(a)$ ($a\in{\bf Z}/N{\bf Z}$).

The mixed pentagon equation \eqref{mixed pentagon}
is a geometric equation introduced by Enriquez \cite{E}.
The series $\varPhi^N_{KZ}$ satisfies the equation, which yields
non-trivial relations among multiple $L$-values.
The generalized double shuffle relation
(the double shuffle relation and the regularization relation)
is a combinatorial relation among multiple $L$-values.
It is formulated as \eqref{double shuffle} for $h=\varPhi^N_{KZ}$.
It is Zhao's  remark  \cite{Z} that for specific $N$'s
the generalized double shuffle relation does not provide all the possible
relations among multiple $L$-values. 

Our main theorem is an implication of 
the generalized double shuffle relation 
from the mixed pentagon equation.

\begin{thm}\label{main}
Let $U\frak F_{N+1}$ be the universal enveloping algebra of 
the free Lie algebra $\frak F_{N+1}$ with variables
$A$ and $B(a)$ ($a\in{\bf Z}/N{\bf Z}$).
Let $h$ be a group-like element in  $U\frak F_{N+1}$
with $c_{B(0)}(h)=0$
satisfying the mixed pentagon equation \eqref{mixed pentagon}
with a group-like series $g\in U\frak F_2$.
Then $h$ also satisfies the generalized double shuffle relation
\eqref{double shuffle}.
\end{thm}

As a consequence we get an embedding from Enriquez' 
cyclotomic associator set $Pseudo(N,\bf Q)$ (definition \ref{Enriquez' set})
defined by the mixed pentagon \eqref{mixed pentagon} and 
the octagon \eqref{octagon} equations into
Racinet's double shuffle set $DMR(N,\bf Q)$ (definition \ref{Racinet's set})
defined by the generalized double shuffle relation \eqref{double shuffle}.

\begin{thm}\label{almost main}
For $N>0$, $Pseudo(N,\bf Q)$ is embedded into $DMR(N,\bf Q)$.
In more detail, we have embeddings of torsors
$$Pseudo_{(a,\mu)}(N,{\bf Q})\subset DMR_{(a,\mu)}(N,\bf Q)$$
for $(a,\mu)\in({\bf Z}/N{\bf Z})^\times \times \bf Q$
and of pro-algebraic groups
$$GRTM_{(\bar 1,1)}(N,{\bf Q})\subset DMR_{(\bar 1,0)}(N,\bf Q)$$
(for $GRTM_{(\bar 1,1)}(N,{\bf Q})$ see definition \ref{GRTM}).
\end{thm}

It might be worthy to note that we do not use the octagon equation
to show $Pseudo(N,{\bf Q})\subset DMR(N,\bf Q)$.
By adding distribution relations \eqref{distribution} to both sides, 
we also get the inclusion
$Psdist_{(a,\mu)}(N,{\bf Q})\subset DMRD_{(a,\mu)}(N,\bf Q)$
(for their definitions see remark \ref{an element in Psdist} and 
\ref{an element in DMRD}).

The motivic fundamental group of the algebraic curve
${\bf G}_m\backslash\mu_N$ is constructed  and discussed in \cite{DG}.
It gives a pro-object of the tannakian category of
mixed Tate motives (constructed in loc.cit.)
over ${\bf Z}[\mu_N,1/N]$,
which yields an action of the motivic Galois group
(: the {\it Galois group} of the tannakian category) 
$Gal^M({\bf Z}[\mu_N,1/N])$ on $\frak F_{N+1}$.
It was  shown that the action is faithful for $N=1$ in \cite{Br}
and $N=2,3,4,8$ in \cite{De08}.
The image of its unipotent part in $Aut \frak F_{N+1}$
is contained in  $GRTMD_{(\bar 1,1)}(N,{\bf Q})$
and $DMRD_{(\bar 1,0)}(N,\bf Q)$,
which follows from a geometric interpretation of the defining equations of $GRTMD_{(\bar 1,1)}(N,{\bf Q})$.
It is a basic problem to ask if they are equal or not.

The contents of the article are as follows:
We recall the mixed pentagon equation in \S \ref{Mixed pentagon equation}
and the generalized double shuffle relation in \S \ref{Double shuffle relation}.
In \S \ref{bar construction} we calculate 
Chen's reduced bar complex
for the Kummer coverings of the moduli spaces $\mathcal M_{0,4}$ and
$\mathcal M_{0,5}$.
Two variable cyclotomic multiple polylogarithms and their associated 
bar elements are introduced in \S \ref{two variables cyclotomic MPL}.
By using them, we prove theorem \ref{main} in \S \ref{proofs of main theorems}.
\S \ref{auxiliary} proves several auxiliary lemmas which are essential
to prove theorem \ref{main}.

\section{Mixed pentagon equation}\label{Mixed pentagon equation}
This section is to recall Enriquez' mixed pentagon equation and
his cyclotomic analogue of the associator set \cite{E}.

Let us fix notations:
For $n\geqslant 2$, the Lie algebra $\frak t_n$
of infinitesimal pure braids is the completed $\bf Q$-Lie algebra
with generators 
$t^{ij}$ ($i\neq j$, $1\leqslant i,j\leqslant n$)
and relations
$$
t^{ij}=t^{ji}, [t^{ij},t^{ik}+t^{jk}]=0 \text{ and }
[t^{ij},t^{kl}]=0  {\text{ for all distinct $i$, $j$, $k$, $l$.}}
$$
We note that $\frak t_2$ is the 1-dimensional abelian Lie algebra
generated by $t^{12}$. 
The element $z_{n}=\sum_{1\leqslant i<j\leqslant n}t^{ij}$
is central in $\frak t_n$.
Put $\frak t^0_{n}$ to be the Lie subalgebra of $\frak t_{n}$
with the same generators except $t^{1n}$.
Then we have $\frak t_{n}=\frak t^0_{n}\oplus{\bf Q}\cdot z_{n}$.
Especially when $n=3$, $\frak t^0_{3}$ is a free Lie algebra $\frak F_{2}$
of rank $2$ with generators $A:=t^{12}$ and $B=t^{23}$.
For a partially defined map
$f:\{1,\dots,m\}\to\{1,\dots,n\}$,
the Lie algebra morphism $\frak t_{n}\to\frak t_{m}:$
$x\mapsto x^f=x^{f^{-1}(1),\dots,f^{-1}(n)}$ is uniquely defined by
$$(t^{ij})^f=\sum_{i'\in f^{-1}(i),j'\in f^{-1}(j)}t^{i'j'}.$$

\begin{defn}[\cite{Dr}]
The associator set $M(\bf Q)$ is defined to be the
set of pairs $(\mu,g)\in {\bf Q}\times\exp\frak F_2=\exp\frak t^0_3$ satisfying
the {\it pentagon equation} in $\exp\frak t^0_4$
\begin{equation}\label{pentagon}
g^{1,2,34}g^{12,3,4}=
g^{2,3,4}g^{1,23,4}g^{1,2,3}.
\end{equation}
and {\it two hexagon equations} in $\exp\frak F_2=\exp\frak t^0_3$
\begin{equation}\label{hexagons}
g(A,B)g(B,A)=1, \
\exp\{\frac{\mu A}{2}\}g(C,A)\exp\{\frac{\mu C}{2}\}g(B,C)\exp\{\frac{\mu B}{2}\}g(A,B)=1
\end{equation}
with $C=-A-B$.
For $\mu\in\bf Q$, the set $M_\mu(\bf Q)$
is the subset of $M(\bf Q)$ with $(\mu,g)\in M(\bf Q)$.
Particularly the set $GRT_1(\bf Q)$ means $M_0(\bf Q)$.

For any field $\bf k$ of characteristic $0$,
$M(\bf k)$ and $GRT(\bf k)$ are defined  in the same way  by replacing $\bf Q$ by $\bf k$. 
\end{defn}

By our notation, the equation \eqref{pentagon} can be read as
$$
g(t^{12},t^{23}+t^{24})g(t^{13}+t^{23},t^{34})=
g(t^{23},t^{34})g(t^{12}+t^{13},t^{24}+t^{34})g(t^{12},t^{23}).
$$
In \cite{Dr} it is shown that
$GRT_1(\bf Q)$ forms a group,
called the  {\it Grothendieck-Teichm\"{u}ller group},
and the associator set $M_\mu(\bf Q)$ with $\mu\in{\bf Q}^\times$
forms a $GRT_1(\bf Q)$-torsor.

\begin{rem}
It is shown in \cite{F10} that the two hexagon equations \eqref{hexagons}
are consequences of the pentagon equation \eqref{pentagon}.
\end{rem}

\begin{rem}
The {\it Drinfel'd associator} $\varPhi_{KZ}=\varPhi_{KZ}(A,B)\in
{\bold C}\langle\langle A,B\rangle\rangle$
is
defined to be the quotient $\varPhi_{KZ}=G_1(z)^{-1}G_0(z)$
where $G_0$ and $G_1$ are the solutions
of the {\it formal KZ (Knizhnik-Zamolodchikov) equation}
$$
\frac{d}{dz}G(z)=\bigl(\frac{A}{z}+\frac{B}{z-1}\bigr)G(z)
$$
such that $G_0(z)\approx z^A$ when $z\to 0$ and
$G_1(z)\approx (1-z)^{B}$ when $z\to 1$  (cf.\cite{Dr}).
The series has the following expression
\begin{align*}
\varPhi_{KZ}(A,B)=1
&+\sum
(-1)^m\zeta(k_1,\cdots,k_m)A^{k_m-1}B\cdots A^{k_1-1}B \\
&+\text{(regularized terms)}
\end{align*}
and the regularized terms are explicitly calculated
to be linear combinations of multiple zeta values
$\zeta(k_1,\cdots,k_m)=L(k_1,\dots,k_m;1,\dots,1)$
in \cite{F03} proposition 3.2.3
by Le-Murakami's method \cite{LM}.
It is shown in \cite{Dr} that the series
belongs to $M_\mu(\bf C)$ with $\mu=2\pi\sqrt{-1}$.
This is achieved by considering monodromy in the real interval $(0,1)$
and the upper half plane in $\mathcal M_{0,4}$ (see \S \ref{bar construction})
and the pentagon formed by
the divisors $y=0$, $x=1$, the exceptional divisor of
the blowing-up at $(1,1)$,
$y=1$ and $x=0$ in ${\mathcal M}_{0,5}$ (see \S\ref{bar construction}).
\end{rem}

For $n\geqslant 2$ and $N\geqslant 1$,
the Lie algebra $\frak t_{n,N}$ is 
the completed $\bf Q$-Lie algebra
with generators 
$$t^{1i} \ (2\leqslant i \leqslant n), \quad
t(a)^{ij} \ (i\neq j, 2\leqslant i,j\leqslant n, a\in\bold Z/N\bold Z)
$$
and relations 
\begin{align*}
&t(a)^{ij}=t(-a)^{ji}, \quad 
[t(a)^{ij},t(a+b)^{ik}+t(b)^{jk}]=0, \\
&[t^{1i}+t^{1j}+\sum_{c\in\bold Z/N\bold Z} t(c)^{ij},t(a)^{ij}]=0, \qquad
[t^{1i},t^{1j}+\sum_{c\in\bold Z/N\bold Z} t(c)^{ij}]=0, \\
&[t^{1i},t(a)^{jk}]=0 \quad
\text{ and }\quad
[t(a)^{ij},t(b)^{kl}]=0
\end{align*}
for all $a$, $b\in\bold Z/N\bold Z$ and
all distinct $i$, $j$, $k$, $l$
($2\leqslant i,j,k,l \leqslant n$).
We note that $\frak t_{n,1}$ is equal to $\frak t_{n}$ for $n\geqslant 2$.
We have a natural injection $\frak t_{n-1,N}\hookrightarrow\frak t_{n,N}$.
The Lie subalgebra $\frak f_{n,N}$ of $\frak t_{n,N}$ 
generated by $t^{1n}$ and
$t(a)^{in}$ ($2\leqslant  i \leqslant n-1$, $a\in\bold Z/N\bold Z$)
is free of rank $(n-2)N+1$
and forms an ideal of $\frak t_{n,N}$.
Actually it shows that $\frak t_{n,N}$ is a semi-direct product of
$\frak f_{n,N}$ and $\frak t_{n-1,N}$.
The element $z_{n,N}=\sum_{1\leqslant i<j\leqslant n}t^{ij}$
with $t^{ij}=\sum_{a\in\bold Z/N\bold Z}t(a)^{ij}$
($2\leqslant i<j\leqslant n$) is central in $\frak t_{n,N}$.
Put $\frak t^0_{n,N}$ to be the Lie subalgebra of $\frak t_{n,N}$
with the same generators except $t^{1n}$.
Then we have $\frak t_{n,N}=\frak t^0_{n,N}\oplus{\bf Q}\cdot z_{n,N}$.
Occasionally we regard $\frak t^0_{n,N}$ as the quotient
$\frak t_{n,N}/{\bf Q}\cdot z_{n,N}$.
Especially when $n=3$, $\frak t^0_{3,N}$ is free Lie algebra $\frak F_{N+1}$
of rank $N+1$ with generators $A:=t^{12}$ and $B(a)=t(a)^{23}$ ($a\in\bold Z/N\bold Z$).

For a partially defined map 
$f:\{1,\dots,m\}\to\{1,\dots,n\}$
such that $f(1)=1$,
the Lie algebra morphism $\frak t_{n,N}\to\frak t_{m,N}$
$x\mapsto x^f=x^{f^{-1}(1),\dots,f^{-1}(n)}$ is uniquely defined by
$$(t(a)^{ij})^f=\sum_{i'\in f^{-1}(i),j'\in f^{-1}(j)}t(a)^{i'j'}
\quad (i\neq j, 2\leqslant i,j\leqslant n)$$
and
\begin{align*}
(t^{1j})^f=&\sum_{j'\in f^{-1}(j)} t^{1j'}+
\frac{1}{2}\sum_{j',j''\in f^{-1}(j)}\sum_{c\in\bold Z/N\bold Z}t(c)^{j'j''} \\
&+\sum_{i'\neq 1\in f^{-1}(1), j'\in f^{-1}(j)}\sum_{c\in\bold Z/N\bold Z}t(c)^{i'j'}
\quad (2\leqslant j\leqslant n).
\end{align*}
Again for a partially defined map 
$g:\{2,\dots,m\}\to\{1,\dots,n\}$,
the Lie algebra morphism $\frak t_{n}\to\frak t_{m,N}$
$x\mapsto x^g=x^{g^{-1}(1),\dots,g^{-1}(n)}$ is uniquely defined by
$$(t^{ij})^g=\sum_{i'\in g^{-1}(i),j'\in g^{-1}(j)}t(0)^{i'j'} \qquad
(i\neq j, 1\leqslant i,j\leqslant n).$$

\begin{defn}[\cite{E}]\label{Enriquez' set}
The cyclotomic associator set $Pseudo(N,\bf Q)$
is defined to be the collection of
$Pseudo_{(a,\mu)}(N,\bf Q)$ for $(a,\mu)\in \bold Z/N{\bold Z}\times \bf Q$
which is defined to be the
set of pairs $(g,h)$
with $g\in \exp\frak F_2$ and
$h=\sum_{W:\text{word}} c_W(h)W\in\exp\frak F_{N+1}$
satisfying $g\in M_\mu(\bf Q)$, $c_{B(0)}(h)=0$,
the {\it mixed pentagon equation} in $\exp\frak t^0_{4,N}$
\begin{equation}\label{mixed pentagon}
h^{1,2,34}h^{12,3,4}=
g^{2,3,4}h^{1,23,4}h^{1,2,3}
\end{equation}
and the {\it octagon equation} in $\exp\frak F_{N+1}$
\begin{align}\label{octagon} 
h\bigl(&A,B(a),B(a+1),\dots,B(a+N-1)\bigr)^{-1}\exp\{\frac{\mu B(a)}{2}\}\cdot
\\ \notag
h&\bigl(C,B(a),B(a-1),\dots,B(a+1-N)\bigr)\exp\{\frac{\mu C}{N}\}\cdot  
\\ \notag
& h\bigl(C,B(0),B(N-1),\dots,B(1)\bigr)^{-1}\exp\{\frac{\mu B(0)}{2}\}\cdot
\\ \notag 
&\quad  h\bigl(A,B(0),B(1),\dots,B(N-1)\bigr)\exp\{\frac{\mu A}{N}\}
=1
\end{align}
with $A+\sum_{a\in\bold Z /N\bold Z}B(a)+C=0$.
\end{defn}

By our notation, each term in the equation \eqref{mixed pentagon} can be read as
\begin{align*}
&h^{1,2,34}=h(t^{12},t^{23}(0)+t^{24}(0),t^{23}(1)+t^{24}(1),\dots,
t^{23}(N-1)+t^{24}(N-1)), \\
&h^{12,3,4}=h(t^{13}+\sum_ct^{23}(c),t^{34}(0),t^{34}(1),\dots,t^{34}(N-1)), \\
&g^{2,3,4}=g(t^{23}(0),t^{34}(0)), \\
&h^{1,23,4}=h(t^{12}+t^{13}+\sum_ct^{23}(c),t^{24}(0)+t^{34}(0),
\dots, t^{24}(N-1)+t^{34}(N-1)),\\
&h^{1,2,3}=h(t^{12},t^{23}(0),t^{23}(1),\dots,t^{23}(N-1)).
\end{align*}

\begin{rem}
In loc.cit. the cyclotomic analogue
$\varPhi^N_{KZ}\in\exp\frak F_{N+1}(\bold C)$
of the Drinfel'd associator
is introduced to be the renormalized holonomy
from $0$ to $1$ of the KZ-like differential equation
$$
\frac{d}{dz}H(z)=\bigl(\frac{A}{z}+\sum_{a\in {\bold Z}/N{\bold Z}}\frac{B(a)}{z-\zeta_N^a}\bigr)H(z)
$$
with $\zeta_N=\exp\{\frac{2\pi\sqrt{-1}}{N}\}$,
i.e.,
$\varPhi^N_{KZ}=H_1^{-1}H_0$ where $H_0$ and $H_1$ are the solutions
such that $H_0(z)\approx z^A$ when $z\to 0$ and
$H_1(z)\approx (1-z)^{B(0)}$ when $z\to 1$  (cf.\cite{E}).
There appear multiple $L$-values \eqref{multiple L-value}
in each of its coefficient;
\begin{align}\label{cyclotomic explicit formula}
\varPhi^N_{KZ}=1+&
\sum (-1)^m L(k_1,\cdots,k_m;
\xi_1,\dots,\xi_m)
A^{k_m-1}B(a_m)\cdots A^{k_1-1}B(a_1) \\ \notag
&+\text{(regularized terms)}
\end{align}
with $\xi_1=\zeta_N^{a_2-a_1}$, \dots,
$\xi_{m-1}=\zeta_N^{a_m-a_{m-1}}$ and
$\xi_m=\zeta_N^{-a_m}$,
where the regularized terms can be explicitly calculated to combinations of
multiple $L$-values by the method of Le-Murakami \cite{LM}.
In \cite{E} it is shown that $(\varPhi_{KZ},\varPhi^N_{KZ})$
belongs to 
$Pseudo_{(-1,2\pi\sqrt{-1})}(N,\bf C)$.
This is achieved by considering monodromy in the pentagon
formed by the divisors $y=0$, $x=1$, the exceptional divisor of
the blowing-up at $(1,1)$,
$y=1$ and $x=0$ in ${\mathcal M}_{0,5}^{(N)}$ (see \S\ref{bar construction})
to get \eqref{mixed pentagon}
and in the octagon formed by $0$, $1$, $\infty$ and $\zeta_N$ in
${\mathcal M}_{0,4}^{(N)}={\bf G}_m\backslash\mu_N$
to get \eqref{octagon}.
\end{rem}

\begin{defn}[\cite{E}]\label{GRTM}
The set $GRTM_{(\bar 1,1)}(N,\bf Q)$ means the subset of
$Pseudo_{(\bar 1,0)}(N,\bf Q)$ satisfying
the {\it special action condition} in $\exp\frak t^0_{3,N}$
\begin{equation}\label{special action condition}
A+
\sum_{a\in\bold Z/N\bold Z}
Ad(\tau_ah^{-1}
)(B(a)) 
+Ad\Bigl(h^{-1} 
\cdot h\bigl(C,B(0),B(N-1),\dots,B(1)\bigr)\Bigr)(C)
=0
\end{equation}
where $\tau_a$ ($a\in\bold Z /N\bold Z $) is the automorphism defined by $A\mapsto A$ and
$B(c)\mapsto B(c+a)$ for all $c\in\bold Z /N\bold Z $.
\end{defn}

In loc.cit. it is shown that $GRTM_{(\bar 1,1)}(N,\bf Q)$
forms a group and and $Pseudo_{(a,\mu)}(N,\bf Q)$ with
$(a,\mu)\in ({\bf Z}/N{\bf Z})^\times \times{\bf Q}^\times$
forms a $GRTM_{(\bar 1,1)}(N,\bf Q)$-torsor.
In the case of $N=1$ we have $g=h$ and then
$M_\mu({\bf Q})=Pseudo_{(1,\mu)}(N,\bf Q)$ and
$GRT_1({\bf Q})=GRTM_{(\bar 1,1)}(N,\bf Q)$.
It is not known for general $N$ whether
$GRTM_{(\bar 1,1)}(N,\bf Q)$ is equal to
$Pseudo_{(\bar 1,0)}(N,\bf Q)$ or not.

Let $N,N'\geqslant 1$ with $N'|N$. Put $d=N/N'$.
The morphism 
$$\pi_{NN'}:\frak t_{n,N}\to\frak t_{n,N'}$$
is defined by
$t^{1i}\mapsto dt^{1i}$ and
$t^{ij}(a)\mapsto t^{ij}(\bar a)$
($i\neq j$, $2\leqslant i,j\leqslant n$, 
$a\in\bold Z/N\bold Z$),
where $\bar a\in\bold Z/N'\bold Z$ means the image of $a$
under the map $\bold Z/N\bold Z\to\bold Z/N'\bold Z$.
The morphism 
$$\delta_{NN'}:\frak t_{n,N}\to\frak t_{n,N'}$$
is defined by
$t^{1i}\mapsto t^{1i}$ and
$t^{ij}(a)\mapsto t^{ij}(a/d)$ if $d|a$ and
$t^{ij}(a)\mapsto 0$ if $d\nmid a$
($i\neq j$, $2\leqslant i,j\leqslant n$, 
$a\in\bold Z/N\bold Z$).
The morphism $\pi_{NN'}$ (resp. $\delta_{NN'}$)
$:\frak t_{n,N}\to\frak t_{n,N'}$ induces the
morphisms ${Pseudo}_{(a,\mu)}(N,{\bf Q})\to
{Pseudo}_{(\bar a,\mu)}(N',{\bf Q})$
which we denote  by the same symbol $\pi_{NN'}$ (resp. $\delta_{NN'}$).
Here $\bar a$ means the image of $a$ by the natural map
$\bold Z/N{\bold Z}\to\bold Z/N'{\bold Z}$.

\begin{rem}\label{an element in Psdist}
In \cite{E}, the {\it distribution relation} in $\exp\frak t^0_{3,N'}$
\begin{equation}\label{distribution}
\pi_{NN'}(h)=\exp\{c_{B(0)}(\pi_{NN'}(h))B(0)\}\delta_{NN'}(h).
\end{equation}
is also discussed and 
${Psdist}_{(a,\mu)}(N,\bf Q)$ (resp.$GRTMD_{(\bar 1,1)}(N,{\bf Q})$)
is defined to be the subset of
${Pseudo}_{(a,\mu)}(N,\bf Q)$ (resp.$GRTM_{(\bar 1,1)}(N,{\bf Q})$)
by adding the distribution relation \eqref{distribution}
in $\exp\frak t^0_{3,N'}$ for all $N'|N$.
In loc.cit.
it is shown that $GRTMD_{(\bar 1,1)}(N,\bf Q)$
forms a group and and $Psdist_{(a,\mu)}(N,\bf Q)$ with
$(a,\mu)\in ({\bf Z}/N{\bf Z})^\times \times{\bf Q}^\times$
forms a $GRTMD_{(\bar 1,1)}(N,\bf Q)$-torsor and
the pair $(\varPhi_{KZ},\varPhi^N_{KZ})$ belongs to it
with $(a,\mu)=(-1,2\pi\sqrt{-1})$.
\end{rem}

\begin{rem}
In \cite{EF} it is proved that the mixed pentagon equation \eqref{mixed pentagon}
implies the distribution relation \eqref{distribution} for $N'=1$
and that the octagon equation \eqref{octagon} follows from the
mixed pentagon equation \eqref{mixed pentagon} 
and the special action condition for $N=2$.
It is also shown that the duality relation shown in \cite{B} is compatible
with the torsor structure of $Psdist(2,\bf Q)$
and a new subtorsor $Psdist^+(2,\bf Q)$ is discussed in \cite{EF}.
\end{rem}

\section{Double shuffle relation}\label{Double shuffle relation}
This section is to recall the generalized double shuffle relation and
Racinet's associated  torsor \cite{R}.

Let us fix notations:
Let $\frak F_{Y_N}$ be the completed graded Lie $\bf Q$-algebra
generated by $Y_{n,a}$ ($n\geqslant 1$ and $a\in\bold Z/N{\bold Z}$)
with $\deg Y_{n,a}=n$.
Put $U\frak F_{Y_N}$ its universal enveloping algebra:
the non-commutative formal series ring with free variables $Y_{n,a}$
($n\geqslant 1$ and $a\in\bold Z/N{\bold Z}$).

Let 
$$\pi_Y: U\frak F_{N+1}\to U\frak F_{Y_N}$$
be the $\bf Q$-linear map between non-commutative formal power series rings
that sends all the words ending in $A$ to zero and the
word $A^{n_m-1}B({a_m})$ $\cdots$$A^{n_1-1}B(a_1)$ 
($n_1,\dots,n_m\geqslant 1$ and $a_1,\dots,a_m\in\bold Z/N{\bold Z}$) to 
$$(-1)^mY_{n_m,-a_m}Y_{n_{m-1},a_m-a_{m-1}}\cdots Y_{n_1,a_2-a_1}.$$
Define the coproduct $\Delta_*$ of $U\frak F_{Y_N}$ by
$$\Delta_* (Y_{n,a})=\sum_{k+l=n,b+c=a} Y_{k,b}\otimes Y_{l,c} \qquad
(n\geqslant 0 \text{ and } a\in\bold Z/N{\bold Z})$$
with $Y_{0,a}:=1$ if $a=0$ and $0$ if $a\neq 0$. 
For $h=\sum_{W:\text{word}} c_W(h) W\in U\frak F_{N+1}$, 
define the series shuffle regularization
$$h_*=h_{\text{corr}}\cdot\pi_Y(h)$$
with the correction term
\begin{equation}\label{correction}
h_{\text{corr}}=\exp\left(\sum_{n=1}^{\infty}
\frac{(-1)^n}{n}c_{A^{n-1}B(0)}(h)Y_{1,0}^n\right).
\end{equation}

\begin{defn}[\cite{R}]\label{Racinet's set}
For $N\geqslant 1$,
the set $DMR(N,\bf Q)$ is defined to be the
set of series $h=\sum_{W:\text{word}} c_W(h) W\in \exp\frak F_{N+1}$ satisfying
$c_A(h)=c_{B(0)}(h)=0$ and 
the {\it generalized double shuffle relation}
\begin{equation}\label{double shuffle}
\Delta_*(h_*)=h_*\widehat\otimes h_*.
\end{equation}
For $(a,\mu)\in{\bf Z}/N{\bold Z}\times\bf Q$, the set
$DMR_{(a,\mu)}(N,\bf Q)$
is defined to be the subset of $DMR(N,\bf Q)$ defined by
\begin{equation}\label{normalization1}
c_{B(ka)}(h)-c_{B(-ka)}(h)=\frac{N-2k}{N-2}\{c_{B(a)}(h)-c_{B(-a)}(h)\}
\end{equation}
for $1\leqslant k \leqslant N/2$ and
\begin{equation}\label{normalization2}
\begin{cases}
c_{AB(0)}(h)=\frac{\mu^2}{24} 
&\text{if  } N=1,2,
\\
c_{B(a)}(h)-c_{B(-a)}(h)=-\frac{(N-2)\mu}{2N} 
&\text{if  } N\geqslant 3.
\end{cases}
\end{equation}
\end{defn}

In loc.cit. it is shown that $DMR_{(\bar 1,0)}(N,\bf Q)$ forms a group
and $DMR_{(a,\mu)}(N,\bf Q)$ with 
$(a,\mu)\in ({\bf Z}/N)^\times\times\bf Q^\times$
forms a $DMR_{(\bar 1,0)}(N,\bf Q)$-torsor.

\begin{rem}\label{an element in DMRD}
In \cite{R},  $DMRD(N,\bf Q)$
(resp. $DMRD_{(a,\mu)}(N,\bf Q)$) is introduced to be the subset of
$DMR(N,\bf Q)$
(resp. $DMR_{(a,\mu)}(N,\bf Q)$) by adding
the distribution relation \eqref{distribution}
in $\exp\frak t^0_{3,N'}=\exp\frak F_{N'+1}$ for all $N'|N$.
The series $\varPhi^N_{KZ}$ belongs to
$DMRD_{(a,\mu)}(N,\bf Q)$
with $(a,\mu)=(-\bar 1,2\pi\sqrt{-1})$ 
because regularized multiple $L$-values satisfy the double shuffle relation
and the distribution relation (loc.cit).
It is shown by Zhao \cite{Z} that for specific $N$'s
all the defining equations of
$DMRD_{(a,\mu)}(N,\bf Q)$ do not provide all the possible relations among
multiple $L$-values.
\end{rem}

\section{Bar constructions}\label{bar construction}
This section reviews the notion of the reduced bar construction
and calculates its $0$-th cohomology for $\mathcal M_{0,4}^{(N)}$ and
$\mathcal M_{0,5}^{(N)}$.

We recall the notion of Chen's reduced bar construction \cite{C}.
Let $(A^\bullet=\oplus_{q=0}^\infty A^q,d)$ be a differential graded algebra (DGA).
The reduced bar complex $\bar B^\bullet(A)$ is the tensor algebra
$\oplus_{r=0}^\infty (\bar A^\bullet)^{\otimes r}$ 
with $\bar A^\bullet=\oplus_{i=0}^\infty \bar A^i$
where $\bar A^0=A^1/dA^0$ and $\bar A^i=A^{i+1}$ ($i>0$).
We denote $a_1\otimes\cdots \otimes a_r$ ($a_i\in {\bar A}^\bullet$) by $[a_1|\cdots|a_r]$.
The degree of elements in $\bar B^\bullet(A)$ is given by the total degree of $\bar A^\bullet$.
Put $Ja=(-1)^{p-1}a$ for $a\in \bar A^p$.
Define
$$d'[a_1|\cdots|a_k]=\sum_{i=1}^k(-1)^i[Ja_1|\cdots|Ja_{i-1}|da_i|
a_{i+1}|\cdots|a_k]$$
and
$$d''[a_1|\cdots|a_k]=\sum_{i=1}^k
(-1)^{i-1}[Ja_1|\cdots|Ja_{i-1}|Ja_i\cdot
a_{i+1}|a_{i+2}|\cdots|a_k].$$
Then $d'+d''$ forms a differential.
The differential and the shuffle product (loc.cit.) give $\bar B^\bullet(A)$
a structure of commutative DGA.
Actually it also forms a Hopf algebra,
whose coproduct $\Delta$ is given by
$$
\Delta([a_1|\cdots|a_r])=\sum_{s=0}^r[a_1|\cdots|a_s]\otimes[a_{s+1}|\cdots|a_r].
$$

For a smooth complex manifold  $\mathcal M$,
$\Omega^\bullet(\mathcal M)$ means the de Rham complex of smooth differential forms on $\mathcal M$
with values in $\bf C$.
We denote the $0$-th cohomology of the reduced bar complex
$\bar B^\bullet(\Omega(\mathcal M))$
with respect to the differential
by $H^0{\bar B}(\mathcal M)$.

Let $\mathcal M_{0,4}$ be the moduli space
$\{(x_1,\cdots,x_4)\in (\bold P_{\bf C}^1)^4|x_i\neq x_j (i\neq j)\}/ PGL_2({\bf C})$
of 4 different points in $\bold P^1$.
It is identified with
$\{z\in\bold P^1_{\bf C} | z\neq 0,1,\infty\}$
by sending $[(0,z,1,\infty)]$ to $z$.
Denote its Kummer $N$-covering 
$${\bf G}_m\backslash\mu_N=\{z\in\bold P^1_{\bf C} | z^N\neq 0,1,\infty\}$$
by $\mathcal M_{0,4}^{(N)}$.
The space $H^0{\bar B}(\mathcal M_{0,4}^{(N)})$
is generated by
$$\omega_0:=d\log(z) \text{ and } \omega_\zeta:=d\log(z-\zeta)
\quad (\zeta\in\mu_N).$$
We have an identification
$H^0{\bar B}(\mathcal M_{0,4}^{(N)})$
with the graded ${\bf C}$-linear dual of
$U\frak F_{N+1}$,
$$
H^0{\bar B}(\mathcal M_{0,4}^{(N)})
\simeq U\frak F_{N+1}^*\otimes{\bf C},
$$
$$\text{by }
\Exp\ \Omega_4^{(N)}
:=\sum X_{i_m}\cdots X_{i_1}\otimes
[\omega_{i_m}|\cdots |\omega_{i_1}]\in
U\frak F_{N+1}\widehat\otimes_{\bf Q} H^0{\bar B}(\mathcal M_{0,4}^{(N)}).$$
Here the sum is taken over $m\geqslant 0$ and $i_1,\cdots,i_m\in\{0\}\cup\mu_N$
and $X_0=A$ and $X_\zeta=B(a)$ when $\zeta=\zeta^a_N$.
It is easy to see that the identification is compatible
with Hopf algebra structures.
We note that the product $l_1\cdot l_2\in H^0{\bar B}(\mathcal M_{0,4}^{(N)})$ 
for $l_1$, $l_2\in H^0{\bar B}(\mathcal M_{0,4}^{(N)})$ is given by
$l_1\cdot l_2(f):=\sum_i l_1(f^{(i)}_1) l_2(f^{(i)}_2)$ for
$f\in U\frak F_{N+1}\otimes{\bf C}$ 
with $\Delta(f)=\sum_i f_1^{(i)}\otimes f_2^{(i)}$.
Occasionally we regard $H^0{\bar B}({\mathcal M}_{0,4}^{(N)})$ as the regular function ring
of $F_{N+1}({\bf C})=\{g\in U\frak F_{N+1}\otimes{\bf C}| g:\text{group-like}\}
=\{g\in U\frak F_{N+1}\otimes{\bf C}| g(0)=1, \Delta(g)=g\otimes g\}$.

Let $\mathcal M_{0,5}$ be the moduli space
$\{(x_1,\cdots,x_5)\in (\bold P_{\bf C}^1)^5|x_i\neq x_j (i\neq j)\}/ PGL_2({\bf C})$
of 5 different points in $\bold P^1$.
It is identified with
$\{(x,y)\in\bold G_m^2 | x\neq 1,y\neq 1,xy\neq 1\}$
by sending $[(0,xy,y,1,\infty)]$ to $(x,y)$.
Denote its Kummer $N^2$-covering
$$\{(x,y)\in\bold G_m^2 | x^N\neq 1,y^N\neq 1,(xy)^N\neq 1\}$$
by $\mathcal M_{0,5}^{(N)}$.
It is identified with $W_{N}/{\bf C}^\times$
by $(x,y)\mapsto(xy,y,1)$
where 
$$W_{N}=\{(z_2,z_3,z_4)\in{\bf G}_m|
z_i^N\neq z_j^N (i\neq j)\}.$$
The space
$H^0{\bar B}(\mathcal M_{0,5}^{(N)})$ is a subspace of the tensor coalgebra
generated by
$$\omega_{1,i}:=d\log{z_i} \text{ and }
\omega_{i,j}(a):=d\log (z_i-\zeta_N^az_j) \quad 
(2\leqslant i,j\leqslant 4,
a\in{\bf Z}/N).$$

\begin{prop}
We have an identification
$$H^0{\bar B}({\mathcal M}_{0,5}^{(N)})\simeq
(U{\bf t}^{0}_{4,N})^*\otimes\bf C.$$
\end{prop}

\begin{pf}
By \cite{K}, $H^0{\bar B}(W_N)$ can be calculated to be the $0$-th cohomology 
$H^0\bar B^\bullet(S)$
of the reduced bar complex of the Orlik-Solomon algebra $S^\bullet$.
The algebra $S^\bullet$ is the (trivial-)differential graded ${\bf C}$-algebra
$S^\bullet=\oplus_{q=0}^\infty S^q$ defined by
generators 
$$\omega_{1,i}=d\log{z_i} \text{ and }
\omega_{i,j}(a)=d\log (z_i-\zeta_N^az_j) \quad (2\leqslant i,j\leqslant 4,
a\in{\bf Z}/N{\bf Z})$$
in degree 1 and
relations
\begin{align*}
&\omega_{i,j}(a)=\omega_{j,i}(-a), \qquad
\omega_{ij}(a)\wedge\{\omega_{ik}(a+b)+\omega_{jk}(b)\}=0, \\
&\{\omega_{1i}+\omega_{1j}+\sum_{c\in\bold Z/N\bold Z}\omega(c)_{ij}\}\wedge\omega(a)_{ij}=0, \\
&\omega_{1i}\wedge\{\omega_{1j}+\sum_{c\in\bold Z/N\bold Z} \omega(c)_{ij}\}=0, \\
&\omega_{1i}\wedge\omega(a)_{jk}=0
\quad \text{ and } \quad
\omega(a)_{ij}\wedge\omega(b)_{kl}=0
\end{align*}
for all $a$, $b\in\bold Z/N\bold Z$ and
all distinct $i$, $j$, $k$, $l$
($2\leqslant i,j,k,l \leqslant n$).
By direct calculation, the element
$$
\sum_{i=2}^4t_{1i}\otimes\omega_{1i}+
\sum_{2\leqslant i<j\leqslant 4, a\in\bold Z/N\bold Z}t_{ij}(a)\otimes\omega_{ij}(a)
\in ({\bf t}_{4,N})^{\deg=1}\otimes S^1
$$
yields a  Hopf  algebra identification of $H^0{\bar B}(W_N)$ with
$(U{\bf t}_{4,N})^*\otimes\bf C$
since both are quadratic.

By the long exact sequence of cohomologies induced from
the ${\bf G}_m$-bundle $W_N\to\mathcal M_{0,5}^{(N)}=W_N/\bf C^\times$,
we get 
$$0\to H^1({\mathcal M}_{0,5}^{(N)})\to H^1(W_N)\to H^1({\bf G}_m)\to 0$$
and  
$$H^i({\mathcal M}_{0,5}^{(N)})\simeq H^i(W_N)\qquad (i\geqslant 2).$$
It yields the identification of the subspace $H^0{\bar B}(\mathcal M_{0,5}^{(N)})$
of $H^0{\bar B}(W_N)$
with $(U{\bf t}^0_{4,N})^*\otimes\bf C$.
\qed
\end{pf}

The above identification is induced from 
$$
\Exp\ \Omega_5^{(N)}:=\sum t_{J_m}\cdots t_{J_1}\otimes
[\omega_{J_m}|\cdots |\omega_{J_1}]\in
U{\bf t}^{0}_{4,N}\widehat\otimes_{\bf Q} H^0{\bar B}({\mathcal M}_{0,5}^{(N)})$$
where the sum is taken over $m\geqslant 0$ and $J_1,\cdots,J_m\in
\{(1,i)|2\leqslant i\leqslant 4\}
\cup\{(i,j,a)|2\leqslant i<j\leqslant 4,a\in{\bf Z}/N{\bf Z}\}$.


Especially the identification between degree 1 terms is given by
\begin{align*}
\Omega_5^{(N)}=\sum_{i=2}^4 t_{1i} &d\log z_i
+\sum_{2\leqslant i<j\leqslant 4}\sum_{a\in{\bf Z}/N{\bf Z}}
t_{i,j}(a) d\log(z_i-\zeta_N^a z_j) \\
&\quad \in \frak t^0_{4,N}\otimes H^1_{DR}(\mathcal M_{0,5}^{(N)}).
\end{align*}
In terms of the coordinate $(x,y)$,
\begin{align*}
\Omega_5^{(N)}
&=t_{12}d\log (xy)+t_{13}d\log y+\sum_a t_{23}(a)d\log y(x-\zeta_N^a)\\
&\qquad\qquad +\sum_a t_{24}(a)d\log (xy-\zeta_N^a)
+\sum_a t_{34}(a)d\log (y-\zeta_N^a)\\
&=t_{12}d\log x+\sum_a t_{23}(a)d\log (x-\zeta_N^a)+(t_{12}+t_{13}+t_{23})d\log y\\
&\qquad\qquad +\sum_a t_{34}(a)d\log (y-\zeta_N^a) 
+\sum_a t_{24}(a)d\log (xy-\zeta_N^a).
\end{align*}
It is easy to see that the identification is compatible
with Hopf algebra structures.
We note again that
the product $l_1\cdot l_2\in H^0{\bar B}(\mathcal M_{0,5}^{(N)})$
for $l_1$, $l_2\in H^0{\bar B}(\mathcal M_{0,5}^{(N)})$ is given by
$l_1\cdot l_2(f):=\sum_i l_1(f^{(i)}_1) l_2(f^{(i)}_2)$ for
$f\in U{\bf t}^{0}_{4,N}\otimes\bf C$ with $\Delta(f)=\sum_i f_1^{(i)}\otimes f_2^{(i)}$
($\Delta$: the coproduct of $U{\bf t}^{0}_{4,N}$).
Occasionally we also regard $H^0{\bar B}({\mathcal M}_{0,5}^{(N)})$ as the regular function ring
of $K_4^N({\bf C})=\{g\in U{\bf t}^{0}_{4,N}\otimes{\bf C}| g:\text{group-like}\}$.

By a generalization of Chen's theory \cite{C} to the case of tangential basepoints,
especially for ${\mathcal M}={\mathcal M_{0,4}^{(N)}}$ or 
${\mathcal M_{0,5}^{(N)}}$,
we have an isomorphism 
$$\rho:H^0{\bar B}(\mathcal M)\simeq I_o(\mathcal M)$$
as algebras over ${\bf C}$ which sends 
$\sum_{I=(i_m,\cdots,i_1)}c_I[\omega_{i_m}|\cdots|\omega_{i_1}]$ 
($c_I\in {\bf C}$) to
$\sum_{I}c_I\text{It}\int_o\omega_{i_m}\circ\cdots\circ\omega_{i_1}$.
Here 
$\sum_{I}c_I\text{It}\int_o\omega_{i_m}\circ\cdots\circ\omega_{i_1}$
means the iterated integral defined by
\begin{equation}\label{iterated integral map}
\sum_{I}c_I\int_{0<t_1< \cdots <t_{m-1}<t_m<1}
\omega_{i_m}({\gamma(t_m)})\cdot
\omega_{i_{m-1}}({\gamma(t_{m-1})})\cdot
\cdots\omega_{i_1}({\gamma(t_1)})
\end{equation}
for all analytic paths $\gamma: (0,1)\to \mathcal M(\bold C)$
starting from the tangential basepoint $o$
(defined by $\frac{d}{dz}$ for $\mathcal M=\mathcal M_{0,4}^{(N)}$ and
defined by $\frac{d}{dx}$ and $\frac{d}{dy}$ for $\mathcal M=\mathcal M_{0,5}^{(N)}$)
at the origin in $\mathcal M$ (for its treatment see also \cite{De}\S 15)
and $I_o(\mathcal M)$
stands for 
the $\bf C$-algebra generated by all such 
homotopy invariant iterated integrals with $m\geqslant 1$ and 
$\omega_{i_1},\dots,\omega_{i_m}\in H^1_{DR}(\mathcal M)$.

\section{Two variable cyclotomic multiple polylogarithms}
\label{two variables cyclotomic MPL}
We introduce cyclotomic multiple polylogarithms,
$Li_{\bold a}(\bar\zeta(z))$ and
$Li_{\bold a,\bold b}(\bar\zeta(x),\bar\eta(y))$,
and their associated  bar elements,
$l_{\bf a}^{\bar\zeta}$ and
$l_{\bf a,\bf b}^{\bar\zeta(x),\bar\eta(y)}$,
which play important roles to prove our main theorems.

For a pair $(\bf a,\bar\zeta)$ with
${\bf a}=(a_1,\cdots,a_k)\in\bold Z^k_{>0}$ and
$\bar{\zeta}=(\zeta_1,\dots,\zeta_k)$
with
$\zeta_i\in\mu_N$: the group of roots of unity in $\bf C$
($1\leqslant i\leqslant k$),
its weight and its depth are defined to be 
$wt({\bf a},\bar\zeta)=a_1+\cdots+a_k$ and
$dp({\bf a},\bar\zeta)=k$ respectively.
Put
$\bar\zeta(x)=(\zeta_1,\dots,\zeta_{k-1},\zeta_k x)$.
Put $z\in\bf C$ with $|z|<1$.
Consider the following complex analytic function,
{\it one variable cyclotomic multiple polylogarithm}
\begin{equation}\label{one variable MPL}
Li_{\bold a}(\bar\zeta(z)):=\underset{0<m_1<\cdots<m_k}{\sum}
\frac{\quad\zeta_1^{m_1}\cdots \zeta_{k-1}^{m_{k-1}}{(\zeta_kz)}^{m_k}}
{m_1^{a_1}\cdots m_{k-1}^{a_{k-1}}m_k^{a_k}}.
\end{equation}
It satisfies the following differential equation
$$
\frac{d}{dz}Li_{\bold a}(\bar\zeta(z))=
\begin{cases}
\frac{1}{z}Li_{(a_1,\cdots,a_{k-1},a_k-1)}(\bar\zeta(z))
&\text{if } a_k\neq 1,\\
\frac{1}{\zeta_k^{-1}-z}Li_{(a_1,\cdots,a_{k-1})}
(\zeta_1,\dots,\zeta_{k-2},\zeta_{k-1}z)
&\text{if } a_k=1,k\neq 1, \\
\frac{1}{\zeta_1^{-1}-z}
&\text{if } a_k=1,k=1. \\
\end{cases}\notag
$$
It gives an iterated integral starting from $o$,
which lies on $I_o(\mathcal M_{0,4}^{(N)})$.
Actually by the map $\rho$ it corresponds to an element of
the $\bf Q$-structure $U\frak F_{N+1}^*$ of
$H^0{\bar B}(\mathcal M_{0,4}^{(N)})$ denoted by $l_{\bf a}^{\bar\zeta}$.
It is expressed as
\begin{equation}\label{Li-expression}
l_{\bf a}^{\bar\zeta}=(-1)^k[\underbrace{\omega_0|\cdots|\omega_0}_{a_k-1}|
\omega_{\zeta^{-1}_k}|
\underbrace{\omega_0|\cdots|\omega_0}_{a_{k-1}-1}|
\omega_{\zeta^{-1}_k\zeta^{-1}_{k-1}}|
\omega_0|\cdots\cdots|\omega_0
|\omega_{\zeta^{-1}_k\cdots\zeta^{-1}_1}].
\end{equation}
By the standard identification $\mu\simeq{\bf Z}/N{\bf Z}$
sending $\zeta_N=\exp\{\frac{2\pi\sqrt{-1}}{N}\}\mapsto 1$,
for a series $\varphi=\sum_{W:\text{word}} c_W(\varphi) W$
it is calculated by
$$l^{\bar\zeta}_{\bf a}(\varphi)=(-1)^k c_{A^{a_k-1}B(-e_k)
A^{a_{k-1}-1}B(-e_k-e_{k-1})\cdots A^{a_1-1}B(-e_k-\cdots-e_1)}(\varphi)
$$
with $\zeta_i=\zeta_N^{e_i}$ ($e_i\in{\bf Z}/N{\bf Z}$).

For ${\bf a}=(a_1,\cdots,a_k)\in\bold Z^k_{>0}$,
${\bf b}=(b_1,\cdots,b_l)\in\bold Z^l_{>0}$,
$\bar{\zeta}=(\zeta_1,\dots,\zeta_k)$,
$\bar{\eta}=(\eta_1,\dots,\eta_l)$
with $\zeta_i,\eta_j\in\mu_N$
and $x,y\in\bf C$ with $|x|<1$ and $|y|<1$,
consider the following complex function,
the {\it two variables multiple polylogarithm}
\begin{equation}\label{two variables MPL}
Li_{\bold a,\bold b}(\bar\zeta(x),\bar\eta(y)):=
\underset{<n_1<\cdots<n_l}{\underset{0<m_1<\cdots<m_k}{\sum}}
\frac{\zeta_1^{m_1}\cdots \zeta_{k-1}^{m_{k-1}}(\zeta_kx)^{m_k}\cdot
\eta_1^{n_1}\cdots\eta_{l-1}^{n_{l-1}}(\eta_ly)^{n_l}}
{m_1^{a_1}\cdots m_{k-1}^{a_{k-1}} m_k^{a_k}\cdot
n_1^{b_1}\cdots n_{l-1}^{b_{l-1}} n_l^{b_l}}.
\end{equation}
It satisfies the following differential equations.

\begin{align}\label{differential equation}
\frac{d}{dx}&Li_{\bold a,\bold b}(\bar\zeta(x),\bar\eta(y))  \\
&=
\begin{cases}
\frac{1}{x}Li_{(a_1,\cdots,a_{k-1},a_k-1),\bold b}(\bar\zeta(x),\bar\eta(y))
\qquad\qquad\quad
\text{if } a_k\neq 1,\\
\frac{1}{\zeta_k^{-1}-x}Li_{(a_1,\cdots,a_{k-1}),\bold b}
(\zeta_1,\dots,\zeta_{k-2},\zeta_{k-1}x,\bar\eta(y))
-\left(\frac{1}{x}+\frac{1}{\zeta_k^{-1}-x}\right)\cdot \\
\quad
Li_{(a_1,\cdots,a_{k-1},b_1),(b_2,\cdots,b_l)}(\zeta_1,\dots\zeta_{k-1},
\zeta_k\eta_1x,\eta_2,\dots,\eta_{l-1},\eta_ly)\\
\qquad\qquad\qquad\qquad\qquad\qquad\qquad\qquad\qquad\quad
\text{if } a_k=1,k\neq 1, l\neq 1,\\
\frac{1}{\zeta_1^{-1}-x}Li_{\bold b}(\eta(y))
-\left(\frac{1}{x}+\frac{1}{\zeta_1^{-1}-x}\right)
Li_{(b_1),(b_2,\cdots,b_l)}(\zeta_1\eta_1x,\eta_2,\dots,\eta_{l-1},\eta_ly) \\
\qquad\qquad\qquad\qquad\qquad\qquad\qquad\qquad\qquad\quad
\text{   if } a_k=1,k= 1, l\neq 1,\\
\frac{1}{\zeta_k^{-1}-x}Li_{(a_1,\cdots,a_{k-1}),b_1}
(\zeta_1,\dots,\zeta_{k-1}x,\eta_1y)
-\left(\frac{1}{x}+\frac{1}{\zeta_k^{-1}-x}\right)\cdot \\
\quad Li_{(a_1,\cdots,a_{k-1},b_1)}(\zeta_1,\dots,\zeta_{k-1},\zeta_k\eta_1xy)
\qquad \ \ 
\text{if } a_k=1,k\neq 1, l=1,\\
\frac{1}{\zeta_1^{-1}-x}Li_{b_1}(\eta_1y)
-\left(\frac{1}{x}+\frac{1}{\zeta_1^{-1}-x}\right)
Li_{b_1}(\zeta_1\eta_1xy)
\
\text{if } a_k=1,k=1, l=1,\\
\end{cases}
\notag 
\end{align}
\begin{align}
\frac{d}{dy}&Li_{\bold a,\bold b}(\bar\zeta(x),\bar\eta(y)) 
\notag  \\
&=
\begin{cases}
\frac{1}{y}Li_{\bold a,(b_1,\cdots,b_{l-1},b_l-1)}(\bar\zeta(x),\bar\eta(y))
&\text{if } b_l\neq 1,\\
\frac{1}{\eta_l^{-1}-y}Li_{\bold a,(b_1,\cdots,b_{l-1})}
(\bar\zeta(x),\eta_1,\dots,\eta_{l-2},\eta_{l-1}y)
&\text{if } b_l=1, l\neq 1,\\
\frac{1}{\eta_1^{-1}-y}Li_{\bold a}(\bar\zeta(\eta_1xy))
&\text{if } b_l=1,l=1. \\
\end{cases}
\notag  
\end{align}
By analytic continuation, the functions 
$Li_{\bf a,\bf b}(\bar\zeta(x),\bar\eta(y))$, 
$Li_{\bf b,\bf a}(\bar\eta(y),\bar\zeta(x))$,
$Li_{\bf a}(\bar\zeta(x))$, $Li_{\bf a}(\bar\zeta(y))$ and 
$Li_{\bf a}(\bar\zeta(xy))$
give iterated integrals starting from $o$,
which lie on $I_o(\mathcal M_{0,5}^{(N)})$.
They correspond to elements of  the
$\bf Q$-structure $(U{\frak t}^0_{4,N})^*$ of
$H^0{\bar B}(\mathcal M_{0,5}^{(N)})$ by the map $\rho$
denoted by
$l_{\bf a,\bf b}^{\bar\zeta(x),\bar\eta(y)}$, 
$l_{\bf b,\bf a}^{\bar\eta(y),\bar\zeta(x)}$,
$l_{\bf a}^{\bar\zeta(x)}$, $l_{\bf a}^{\bar\eta(y)}$  and $l_{\bf a}^{\bar\zeta(xy)}$ respectively.
Note that they are expressed as
\begin{equation}\label{expression}
\sum_{I=(i_m,\cdots,i_1)}c_I[\omega_{i_m}|\cdots|\omega_{i_1}]
\end{equation}
for some $m\in\bf N$
with $c_I\in \bold Q$ and $\omega_{i_j}\in \{\frac{dx}{x},\frac{dx}{\zeta-x},\frac{dy}{y},\frac{dy}{\zeta-y},
\frac{xdy+ydx}{\zeta-xy} (\zeta\in\mu_N)\}$.

\section{Proofs of main theorems}\label{proofs of main theorems}
This section gives proofs of theorem \ref{main} and theorem \ref{almost main}.

{\bf Proof of theorem \ref{main}.}
Let ${\bf a}=(a_1,\dots,a_k)\in\bold Z^k_{>0}$,
${\bf b}=(b_1,\dots,b_l)\in\bold Z^l_{>0}$,
$\bar{\zeta}=(\zeta_1,\dots,\zeta_k)$ and
$\bar\eta=(\eta_1,\dots,\eta_l)$
with
$\zeta_i,\eta_j\in\mu_N\subset \bf C$ ($1\leqslant i\leqslant k$
and $1\leqslant j\leqslant l$). 
Put
$\bar\zeta(x)=(\zeta_1,\dots,\zeta_{k-1},\zeta_k x)$ and
$\bar\eta(y)=(\eta_1,\dots,\eta_{l-1},\eta_ly)$.
Recall that multiple polylogarithms satisfy the following analytic identity,
the series shuffle formula in
$I_o(\mathcal M_{0,5}^{(N)})$:
$$
Li_{\bold a}(\bar\zeta(x))\cdot Li_{\bold b}(\bar\eta(y))
={\underset{\sigma\in Sh^{\leqslant}(k,l)}{\sum}}
Li_{\sigma(\bold a,\bold b)}^{\sigma(\bar\zeta(x),\bar\eta(y))}.
$$
Here
$
Sh^{\leqslant}(k,l):={\cup}^\infty_{N=1}\{
\sigma:\{1,\cdots,k+l\}\to\{1,\cdots,N \}| \sigma{\text{ is onto}},
\sigma(1)<\cdots<\sigma(k), \sigma(k+1)<\cdots<\sigma(k+l)
\},
$
$\sigma(\bold a,\bold b):=(c_1,\cdots,c_{N})$ with
$$
c_i=
\begin{cases}
a_s+b_{t-k}  &\text{if } \sigma^{-1}(i)=\{s,t\} \text{ with } s<t , \\
a_s     &\text{if } \sigma^{-1}(i)=\{s\}  \quad \text{with } s\leqslant k,\\
b_{s-k} &\text{if } \sigma^{-1}(i)=\{s\}  \quad \text{with } s> k,\\
\end{cases}
$$
and 
$\sigma(\bar\zeta(x),\bar\eta(y)):=(z_1,\dots,z_N)$ with
$$
z_i=
\begin{cases}
x_sy_{t-k}  &\text{if } \sigma^{-1}(i)=\{s,t\} \text{ with } s<t , \\
x_s     &\text{if } \sigma^{-1}(i)=\{s\}  \quad \text{with } s\leqslant k,\\
y_{s-k} &\text{if } \sigma^{-1}(i)=\{s\}  \quad \text{with } s> k,\\
\end{cases}
$$
for $x_i=\zeta_i$ ($i\neq k$), $\zeta_k x$ ($i=k$) and
$y_j=\eta_j$ ($j\neq l$), $\eta_j y$ ($j=l$).
Since $\rho$ is an embedding of algebras,
the above analytic identity immediately implies the algebraic identity,
the series shuffle formula in the $\bf Q$-structure $(U{\frak t}^0_{4,N})^*$
of $H^0{\bar B}(\mathcal M_{0,5}^{(N)})$
\begin{equation}\label{series shuffle for bar}
l_{\bold a}^{\bar\zeta(x)}\cdot l_{\bold b}^{\bar\eta(y)}
={\underset{\sigma\in Sh^{\leqslant}(k,l)}{\sum}}
l_{\sigma(\bold a,\bold b)}^{\sigma(\bar\zeta(x),\bar\eta(y))}.
\end{equation}

Let $(g,h)$ be a pair in theorem \ref{main}.
By the group-likeness of $h$, i.e. $h\in\exp\frak F_{N+1}$,
the product $h^{1,23,4}h^{1,2,3}$ is group-like,
i.e. belongs to $\exp\frak t^0_{4,N}$.
Hence $\Delta{(h^{1,23,4}h^{1,2,3})}=
(h^{1,23,4}h^{1,2,3})\widehat{\otimes} (h^{1,23,4}h^{1,2,3})$,
where $\Delta$ is the standard coproduct of $U\frak t^0_{4,N}$.
Therefore
\begin{align*}
l_{\bold a}^{\bar\zeta(x)}\cdot l_{\bold b}^{\bar\eta(y)}
(h^{1,23,4}h^{1,2,3})=&
(l_{\bold a}^{\bar\zeta(x)}\widehat{\otimes}l_{\bold b}^{\bar\eta(y)})
(\Delta(h^{1,23,4}h^{1,2,3}))\\
=&
l_{\bold a}^{\bar\zeta(x)}(h^{1,23,4}h^{1,2,3})\cdot 
l_{\bold b}^{\bar\eta(y)}(h^{1,23,4}h^{1,2,3}).
\end{align*}
Evaluation of the equation \eqref{series shuffle for bar}
at the group-like element 
$h^{1,23,4}h^{1,2,3}$
gives the series shuffle formula
\begin{equation}\label{series shuffle for l}
l_{\bold a}^{\bar\zeta}(h)\cdot l_{\bold b}^{\bar\eta}(h)
={\underset{\sigma\in Sh^{\leqslant}(k,l)}{\sum}}
l_{\sigma(\bold a,\bold b)}^{\sigma(\bar\zeta,\bar\eta)}(h)
\end{equation}
for admissible pairs
\footnote{
A pair $({\bf a},\bar\zeta)$ with 
${\bf a}=(a_1,\cdots,a_k)$ and $\bar\zeta=(\zeta_1,\dots,\zeta_k)$
is called {\it admissible} if $(a_k,\zeta_k)\neq(1,1)$.
}
$({\bf a},\bar\zeta)$ and $({\bf b},\bar\eta)$
by lemma \ref{lemma3} and lemma \ref{lemma4} below
because the group-likeness and \eqref{mixed pentagon} for $h$
implies $c_0(h)=1$ and $c_A(h)=0$.

By putting $l^{1,S}_1(h):=-T$ and
$l^{\bar\zeta,S}_{\bf a}(h):=l^{\bar\zeta}_{\bf a}(h)$ for all admissible 
pairs $({\bf a},\bar\zeta)$,
the series regularized value $l^{\bar\zeta,S}_{\bf a}(h)$ 
in ${\bf Q}[T]$ ($T$: a parameter which stands for $\log z$. cf. \cite{R})
for a non-admissible pair $({\bf a},\bar\zeta)$
is uniquely determined
in such a way (cf.\cite{AK}) that the above series shuffle formulae 
remain valid for
$l^{\bar\zeta,S}_{\bf a}(h)$ 
with all pairs $({\bf a},\bar\zeta)$.

Define the integral regularized value $l^{\bar\zeta,I}_{\bf a}(h)$
in ${\bf Q}[T]$
for all pairs $({\bf a},\bar\zeta)$ 
by $l^{\bar\zeta,I}_{\bf a}(h)=\l^{\bar\zeta}_{\bf a}(e^{TB(0)}h)$.
Equivalently $l^{\bar\zeta,I}_{\bf a}(h)$ for any pair $({\bf a},\bar\zeta)$
can be uniquely defined in such a way that
the iterated integral shuffle formulae (loc.cit)
remain valid 
for all  pairs $({\bf a},\bar\zeta)$
with $l^{1,I}_1(h):=-T$ and
$l^{\bar\zeta,I}_{\bf a}(h):=l^{\bar\zeta}_{\bf a}(h)$ for all admissible 
pairs $({\bf a},\bar\zeta)$
because they hold for admissible pairs by the group-likeness of $h$
(cf. loc.cit).

Let ${\mathbb L}$ be the $\bf Q$-linear map from
${\bf Q}[T]$ to itself defined via the generating function:
\begin{align}\label{generating function}
{\mathbb L}(\exp Tu)
&=\sum_{n=0}^{\infty}{\mathbb L}(T^n)\frac{u^n}{n!}
=\mbox{exp}\left\{-\sum_{n=1}^{\infty}{l^{1,I}_n(h)}\frac{u^n}{n}\right\} \\
&\left(=\mbox{exp}\left\{Tu-\sum_{n=1}^{\infty}{l^1_n(h)}\frac{u^n}{n}\right\}\right).
\notag
\end{align}

\begin{prop}\label{regularization relation}
Let $h$ be an element as in theorem \ref{main}.
Then the regularization relation holds, i.e.
$\l^{\bar\zeta,S}_{\bf a}(h)={\mathbb L}\bigl(l^{\bar\zeta,I}_{\bf a}(h)\bigr)$
for all pairs $({\bf a},\bar\zeta)$.
\end{prop}

\begin{pf}
We may assume that $({\bf a},\bar\zeta)$ is non-admissible
because the proposition is trivial if it is admissible.
Put $1^n=(\underbrace{1,1,\cdots,1}_{n})$.
When ${\bf a}=1^n$ and ${\bar\zeta}=\bar{1}^n$,
the proof is given by the same argument to \cite{F08} as follows:
By the series shuffle formulae,
$$
\sum_{k=0}^m (-1)^kl^{\bar 1,S}_{k+1}(h)\cdot l^{\bar 1^{m-k},S}_{1^{m-k}}(h)
=(m+1)l^{\bar 1^{m+1},S}_{1^{m+1}}(h)
$$
for $m\geqslant 0$. Here we put $l^{\emptyset,S}_\emptyset(h)=1$.
This means
$$
\sum_{k,l\geqslant 0} (-1)^kl^{\bar 1,S}_{k+1}(h)\cdot 
l^{\bar 1^l,S}_{1^{l}}(h)u^{k+l}
=\sum_{m\geqslant 0}(m+1)l^{\bar 1^{m+1},S}_{1^{m+1}}(h)u^m.
$$
Put
$f(u)=\sum_{n\geqslant 0}l^{\bar 1^n,S}_{1^{n}}(h)u^{n}$.
Then the above equality can be read as
$$
\sum_{k\geqslant 0} (-1)^kl^{\bar 1,S}_{k+1}(h)u^k=
\frac{d}{du}\log f(u).
$$
Integrating and adjusting constant terms gives
$$
\sum_{n\geqslant 0}l^{\bar 1^n,S}_{1^n}(h)u^{n}
=\exp\left\{
-\sum_{n\geqslant 1} (-1)^nl^{\bar 1,S}_{n}(h)\frac{u^n}{n}
\right\}
=\exp\left\{
-\sum_{n\geqslant 1} (-1)^nl^{\bar 1,I}_{n}(h)\frac{u^n}{n}
\right\}
$$
because $l^{\bar 1,S}_n(h)=l^{\bar 1,I}_n(h)=l^1_n(h)$ for $n>1$ and $l^{\bar 1,S}_1(h)=l^{\bar 1,I}_1(h)=-T$.
Since $l^{\bar 1^m,I}_{1^m}(h)=\frac{(-T)^m}{m!}$,
we get $\l^{\bar 1^m,S}_{1^m}(h)={\mathbb L}\bigl(l^{\bar 1^m,I}_{1^m}(h)\bigr)$.

When $({\bf a},\bar\zeta)$
is of the form $({\bf a'}1^l,\bar{\zeta'}\bar{1^l})$
with $({\bf a'},\bar{\zeta'})$ admissible,
the proof is given by the following induction on $l$.
By \eqref{series shuffle for bar},
$$
l_{\bold a'}^{\bar{\zeta'}(x)}(h')\cdot 
l_{1^l}^{\bar{1^l}(y)}(h')
={\underset{\sigma\in Sh^{\leqslant}(k,l)}{\sum}}
l_{\sigma(\bold a',1^l)}^{\sigma({\bar{\zeta'}(x)},{\bar{1^l}(y)})}(h')
$$
for $h'=e^{T\{t^{23}(0)+t^{24}(0)+t^{34}(0)\}}h^{1,23,4}h^{1,2,3}$
with $k=dp({\bf a'})$.
The group-likeness and \eqref{mixed pentagon} for $h$
implies $c_0(h)=1$ and $c_A(h)=0$ and
the group-likeness and our assumption $c_{B(0)}(h)=0$ implies
$c_{B(0)^n}(h)=0$ for $n\in\bold Z_{>0}$.
Hence by lemma \ref{lemma5} and lemma \ref{lemma6},
$$
l_{\bold a'}^{\bar{\zeta'}}(h)\cdot 
l_{1^l}^{\bar{1^l},I}(h)
={\underset{\sigma\in Sh^{\leqslant}(k,l)}{\sum}}
l_{\sigma(\bold a',1^l)}^{\sigma({\bar{\zeta'}},{\bar{1^l}}),I}(h).
$$
Then by our induction assumption, taking the image by the map $\mathbb L$ gives 
$$
l_{\bold a'}^{\bar{\zeta'}}(h)\cdot 
l_{1^l}^{\bar{1^l},S}(h)
={\mathbb L}\bigl(l^{\bar{\zeta'}{\bar{1^l}},I}_{{\bf a'}1^l}(h)\bigr)
+{\underset{\sigma\neq id\in Sh^{\leqslant}(k,l)}{\sum}}
l_{\sigma(\bold a',1^l)}^{\sigma({\bar{\zeta'}},{\bar{1^l}}),S}(h).
$$
Since $l^{{\bar{\zeta'}},S}_{\bf a'}(h)$ and 
$l_{1^l}^{\bar{1^l},S}(h)$ satisfy the series shuffle formula,
${\mathbb L}\bigl(l^{{\bar{\zeta}},I}_{\bf a}(h)\bigr)$ must be equal to
$\l^{\bar\zeta,S}_{\bf a}(h)$,
which concludes proposition \ref{regularization relation}.
\qed
\end{pf}

Embed $U\frak F_{Y_{N}}$ 
into $U\frak F_{N+1}$
by sending $Y_{m,a}$ to $-A^{m-1}B(-a)$.
Then by the above proposition,
\begin{align*}
\l^{\bar\zeta,S}_{\bf a}(h)
&
={\mathbb L}(l^{\bar\zeta,I}_{\bf a}(h))
={\mathbb L}(l^{\bar\zeta}_{\bf a}(e^{TB(0)}h))
=l^{\bar\zeta}_{\bf a}\left({\mathbb L}(e^{TB(0)}\pi_Y(h))\right) \\
&=l^{\bar\zeta}_{\bf a}(\mbox{exp}
\left\{-\sum_{n=1}^{\infty}{l^{1,I}_n(h)}
\frac{B(0)^n}{n}\right\}\cdot\pi_Y(h)) \\
&=l^{\bar\zeta}_{\bf a}(\mbox{exp}
\left\{-TY_{1,0}+\sum_{n=1}^{\infty}\frac{(-1)^n}{n}
c_{A^{n-1}B(0)}(h)Y^n_{1,0}\right\}
\cdot\pi_Y(h)) \\
&=l^{\bar\zeta}_{\bf a}(e^{-TY_{1,0}}h_*)
\end{align*}
for all $(\bf a,\bar\zeta)$ because $l^1_1(h)=0$.
As for the third equality we use
$$({\mathbb L}\otimes_{\bf Q} id)\circ (id\otimes_{\bf Q} l^{\bar\zeta}_{\bf a})=
(id\otimes_{\bf Q} l^{\bar\zeta}_{\bf a})\circ ({\mathbb L}\otimes_{\bf Q} id)
\text{ on } {\bf Q}[T]\otimes_{\bf Q} U\frak F_{N+1}.$$
All $\l^{\bar\zeta,S}_{\bf a}(h)$'s satisfy 
the series shuffle formulae \eqref{series shuffle for l},
so the $l^{\bar\zeta}_{\bf a}(e^{-TY_{1,0}}h_*)$'s do also.
By putting $T=0$, we get that $l^{\bar\zeta}_{\bf a}(h_*)$'s also satisfy the series shuffle formulae for all $\bf a$.
Therefore $\Delta_*(h_*)=h_*\widehat\otimes h_*$.
This completes the proof of theorem \ref{main}.
\qed

{\bf Proof of theorem \ref{almost main}.}
The first statement follows from theorem \ref{main}.

Let $(g,h)\in Pseudo_{(a,\mu)}(N,{\bf Q})$
with $(a,\mu)\in ({\bf Z}/N{\bf Z})^\times\times{\bf Q}$.
By comparing the coefficient of $B(a)$ in 
the octagon equation \eqref{octagon},
$$-c_{B(0)}(h)+\frac{\mu}{2}-c_A(h)+c_{B(0)}(h)-\frac{\mu}{N}
+c_A(h)-c_{B(-a)}(h)+c_{B(a)}(h)=0.
$$
Thus
$c_{B(a)}(h)-c_{B(-a)}(h)=(\frac{1}{N}-\frac{1}{2})\mu$.

Next by comparing  the coefficient of $B(ka)$ in \eqref{octagon}
for $2\leqslant k\leqslant N/2$,
$$-c_{B((k-1)a)}(h)-c_A(h)+c_{B(-(k-1)a)}(h)-\frac{\mu}{N}
+c_A(h)-c_{B(-ka)}(h)+c_{B(ka)}(h)=0.
$$
Thus $c_{B(ka)}(h)-c_{B(-ka)}(h)=c_{B((k-1)a)}(h)-c_{B(-(k-1)a)}(h)
+\frac{\mu}{N}$.

By combining these equations we get
\eqref{normalization1} and \eqref{normalization2}
for $N\geqslant 3$.
Since we have $c_{AB}(g)=\frac{\mu^2}{24}$ for $g\in M_\mu(\bf Q)$,
we have \eqref{normalization2}
for $N=1,2$ by $c_{AB}(g)=c_{AB(0)}(h)$.
\qed

\section{Auxiliary lemmas}\label{auxiliary}
We prove all lemmas which are required to prove theorem \ref{main}
in the previous section.  

\begin{lem}\label{lemma3}
Let $h\in U \frak F_{N+1}$ with $c_0(h)=1$
\footnote{
The symbol $c_0(h)$ stands for the constant term of $h$.
}
and $c_{A}(h)=0$.
Then
$$l^{\bar\zeta(x)}_{\bf a}(h^{1,23,4}h^{1,2,3})=l^{\bar\zeta}_{\bf a}(h),$$
$$l^{\bar\zeta(y)}_{\bf a}(h^{1,23,4}h^{1,2,3})=l^{\bar\zeta}_{\bf a}(h),$$
$$l^{\bar\zeta(xy)}_{\bf a}(h^{1,23,4}h^{1,2,3})=l^{\bar\zeta}_{\bf a}(h),$$
$$l^{\bar\zeta(x),\bar\eta(y)}_{\bf a,\bf b}(h^{1,23,4}h^{1,2,3})
=l^{\bar\zeta \bar\eta}_{\bf a\bf b}(h)$$
for any pairs $(\bf a,\bar\zeta)$ and $(\bf b,\bar\eta)$.
\end{lem}

\begin{pf}
Put $U\frak t^0_{4,N}$ the universal enveloping algebra of $\frak t^0_{4,N}$.
Consider the map 
$\mathcal M^{(N)}_{0,5}\to\mathcal M^{(N)}_{0,4}$ induced from
$\mathcal M_{0,5}\to\mathcal M_{0,4}:
[(x_1,\cdots,x_5)]\mapsto [(x_1,x_2,x_3,x_5)]$.
This yields the projection $p_{4}: U\frak t^0_{4,N} \twoheadrightarrow U\frak F_{N+1}$
sending  $t^{14},t^{24}(a),t^{34}(a)\mapsto 0$,
$t^{12}\mapsto A$ and $t^{23}(a)\mapsto B(a)$
($a\in{\bf Z}/N{\bf Z}$).
Express $l^{\bar\zeta}_{\bf a}$ as \eqref{Li-expression}.
Since $(p_4\otimes id)(\Exp\Omega^{(N)}_5)=\Exp\Omega^{(N)}_4(x)\in
U\frak F_{N+1}\widehat\otimes_{\bf Q} H^0{\bar B}({\mathcal M}^{(N)}_{0,5})
\simeq H^0{\bar B}({\mathcal M}^{(N)}_{0,4})^*\widehat\otimes_{\bf C} H^0{\bar B}({\mathcal M}^{(N)}_{0,5})$,
it induces the map
$$p^*_4:H^0{\bar B}({\mathcal M}^{(N)}_{0,4})\to H^0{\bar B}({\mathcal M}^{(N)}_{0,5})$$
which gives
$p^*_4([\frac{dz}{z}])=[\frac{dx}{x}]$ and
$p^*_4([\frac{dz}{\zeta^a_N-z}])=[\frac{dx}{\zeta^a_N-x}]$.
Hence 
$$p_4^*(l^{\bar\zeta}_{\bf a})=l^{\bar\zeta(x)}_{\bf a}.$$
Then $l^{\bar\zeta(x)}_{\bf a}(h^{1,23,4}h^{1,2,3})=
l^{\bar\zeta}_{\bf a}(p_4(h^{1,23,4}h^{1,2,3}))
=l^{\bar\zeta}_{\bf a}(h)$
because $p_4(h^{1,23,4})=0$ by our assumption  $c_A(h)=0$.

Next consider the map 
$\mathcal M^{(N)}_{0,5}\to\mathcal M^{(N)}_{0,4}$ induced from
$\mathcal M_{0,5}\to\mathcal M_{0,4}$: 
$[(x_1,\cdots,x_5)]\mapsto [(x_1,x_3,x_4,x_5)]$.
This induces the projection $p_{2}: U\frak t^0_{4,N} \twoheadrightarrow U\frak F_{N+1}$
sending $t^{12},t^{23}(a), t^{24}(a)\mapsto 0$, 
$t^{12}+t^{13}+t^{23}\mapsto A$ and $t^{34}(a)\mapsto B(a)$
($a\in{\bf Z}/N{\bf Z}$).
Since $(p_2\otimes id)(\Exp\Omega^{(N)}_5)=\Exp\Omega^{(N)}_4(y)\in
U\frak F_{N+1}\widehat\otimes_{\bf Q} H^0{\bar B}({\mathcal M}^{(N)}_{0,5})
\simeq H^0{\bar B}({\mathcal M}^{(N)}_{0,4})^*{\widehat\otimes_{\bf C}}$
$H^0{\bar B}({\mathcal M}^{(N)}_{0,5})$,
it induces the map
$$p^*_2:H^0{\bar B}({\mathcal M}^{(N)}_{0,4})\to H^0{\bar B}({\mathcal M}^{(N)}_{0,5})$$
which gives
$p^*_2([\frac{dz}{z}])=[\frac{dy}{y}]$ and
$p^*_2([\frac{dz}{\zeta^a_N-z}])=[\frac{dy}{\zeta^a_N-y}]$.
Hence 
$$p_2^*(l^{\bar\zeta}_{\bf a})=l^{\bar\zeta(y)}_{\bf a}.$$
Then $l^{\bar\zeta(y)}_{\bf a}(h^{1,23,4}h^{1,2,3})=
l^{\bar\zeta}_{\bf a}(p_2(h^{1,23,4}h^{1,2,3}))
=l^{\bar\zeta}_{\bf a}(h)$
because $p_2(h^{1,2,3})=0$.

Similarly consider the map 
$\mathcal M^{(N)}_{0,5}\to\mathcal M^{(N)}_{0,4}$ induced from
$\mathcal M_{0,5}\to\mathcal M_{0,4}$:
$[(x_1,\cdots,x_5)]\mapsto [(x_1,x_2,x_4,x_5)]$.
This induces the projection $p_{3}: U\frak t^0_{4,N} \twoheadrightarrow U\frak F_{N+1}$
sending  $t^{13},t^{23}(a),t^{34}(a)\mapsto 0$,
$t^{12}\mapsto A$ and $t^{24}(a)\mapsto B(a)$
($a\in{\bf Z}/N{\bf Z}$).
Since $(p_3\otimes id)(\Exp\Omega^{(N)}_5)=\Exp\Omega^{(N)}_4(xy)\in
U\frak F_{N+1}\widehat\otimes_{\bf Q} H^0{\bar B}({\mathcal M}^{(N)}_{0,5})
\simeq H^0{\bar B}({\mathcal M}^{(N)}_{0,4})^*\widehat\otimes_{\bf C} H^0{\bar B}({\mathcal M}^{(N)}_{0,5})$,
it induces the map
$$p^*_3:H^0{\bar B}({\mathcal M}^{(N)}_{0,4})\to H^0{\bar B}({\mathcal M}^{(N)}_{0,5})$$
which gives
$p^*_3([\frac{dz}{z}])=[\frac{dx}{x}+\frac{dy}{y}]$ and
$p^*_3([\frac{dz}{\zeta^a_N-z}])=[\frac{xdy+ydx}{\zeta^a_N-xy}]$.
Hence 
$$p_3^*(l^{\bar\zeta}_{\bf a})=l^{\bar\zeta(xy)}_{\bf a}.$$
Then $l^{\bar\zeta(xy)}_{\bf a}(h^{1,23,4}h^{1,2,3})=
l^{\bar\zeta}_{\bf a}(p_3(h^{1,23,4}h^{1,2,3}))
=l^{\bar\zeta}_{\bf a}(h)$
because $p_3(h^{1,2,3})=0$ by our assumption $c_A(h)=0$.

Consider the embedding of Hopf algebras
$i_{1,2,3}:U\frak F_{N+1}\hookrightarrow U\frak t^0_{4,N}$
sending $A\mapsto t^{12}$ and $B(a)\mapsto t^{23}(a)$
along the divisor $\{y=0\}$.
Since
$(i_{1,2,3}\otimes id)(\Exp\Omega^{(N)}_4)=\Exp\Omega^{(N)}_4(z)^{1,2,3}
\in U\frak t^0_{4,N}\widehat\otimes_{\bf Q} H^0{\bar B}({\mathcal M}^{(N)}_{0,4})
\simeq H^0{\bar B}({\mathcal M}^{(N)}_{0,5})^*\widehat\otimes_{\bf C} H^0{\bar B}({\mathcal M}^{(N)}_{0,4})$,
it induces the map
$$i^*_{1,2,3}:H^0{\bar B}({\mathcal M}^{(N)}_{0,5})\to H^0{\bar B}({\mathcal M}^{(N)}_{0,4})$$ 
which gives $i^*_{1,2,3}([\frac{dy}{y}])=i^*_{1,2,3}([\frac{dy}{\zeta^a_N-y}])
=i^*_{1,2,3}([\frac{xdy+ydx}{\zeta^a_N-xy}])=0$.
Express $l^{\bar\zeta(x),\bar\eta(y)}_{\bf a,\bf b}$ and $l^{\bar\zeta(xy)}_{\bf a}$ as \eqref{expression}.
In the expression each term contains at least one
$\frac{dy}{y}$, $\frac{dy}{\zeta^a_N-y}$ or $\frac{xdy+ydx}{\zeta^a_N-xy}$.
Therefore we have 
$$i_{1,2,3}^*(l^{\bar\zeta(x),\bar\eta(y)}_{\bf a,\bf b})=0 \text{  and  }
i_{1,2,3}^*(l^{\bar\zeta(xy)}_{\bf a})=0.$$
Thus $l^{\bar\zeta(x),\bar\eta(y)}_{\bf a,\bf b}(h^{1,2,3})=
i_{1,2,3}^*(l^{\bar\zeta(x),\bar\eta(y)}_{\bf a,\bf b})(h)=0$
and $l^{\bar\zeta(xy)}_{\bf a}(h^{1,2,3})=
i_{1,2,3}^*(l^{\bar\zeta(xy)}_{\bf a})(h)=0$.

Next consider the embedding of Hopf algebras
$i_{1,23,4}:U\frak F_{N+1}\hookrightarrow U\frak t^0_{4,N}$
sending $A\mapsto t^{12}+t^{13}+t^{23}$ and $B(a)\mapsto t^{24}(a)+t^{34}(a)$
(geometrically caused by the divisor $\{x=1\}$.)
Since
$(i_{1,23,4}\otimes id)(\Exp\Omega^{(N)}_4)=\Exp\Omega^{(N)}_4(z)^{1,23,4}
\in U\frak t^0_{4,N}\widehat\otimes_{\bf Q} H^0{\bar B}({\mathcal M}^{(N)}_{0,4})
\simeq H^0{\bar B}({\mathcal M}^{(N)}_{0,5})^*\widehat\otimes_{\bf C} H^0{\bar B}({\mathcal M}^{(N)}_{0,4})$,
it induces the map
$$i^*_{1,23,4}:H^0{\bar B}({\mathcal M}^{(N)}_{0,5})\to H^0{\bar B}({\mathcal M}^{(N)}_{0,4})$$
which gives
$i^*_{1,23,4}([\frac{dx}{x}])=0$,
$i^*_{1,23,4}([\frac{dx}{\zeta^a_N-x}])=[\frac{dz}{\zeta^a_N-z}]$,
$i^*_{1,23,4}([\frac{dy}{y}])=[\frac{dz}{z}]$,
$i^*_{1,23,4}([\frac{dy}{\zeta^a_N-y}])=[\frac{dz}{\zeta^a_N-z}]$ and
$i^*_{1,23,4}([\frac{xdy+ydx}{\zeta^a_N-xy}])=[\frac{dz}{\zeta^a_N-z}]$.
As is same to the proof of \cite{F08} lemma 5.1,
$$i_{1,23,4}^*(l^{\bar\zeta(x),\bar\eta(y)}_{\bf a,\bf b})
=l^{\bar\zeta\bar\eta}_{\bf a\bf b} \text{  and  }
i_{1,23,4}^*(l^{\bar\zeta(xy)}_{\bf a})=l^{\bar\zeta}_{\bf a}$$
can be deduced by induction on weight.
Thus $l^{\bar\zeta(x),\bar\eta(y)}_{\bf a,\bf b}(h^{1,23,4})=l^{\bar\zeta\bar\eta}_{\bf a\bf b}(h)$.
Let $\delta$ be the coproduct of $H^0{\bar B}({\mathcal M}^{(N)}_{0,5})$.
Express $\delta(l^{\bar\zeta(x),\bar\eta(y)}_{\bf a,\bf b})=\sum_i l'_i\otimes l''_i$
with $\deg l'_i={m'_i}$ and $\deg l''_i={m''_i}$ for some $m'_i$ and $m''_i$
such that $m'_i+m''_i=wt({\bf a},\bar\zeta)+wt({\bf b},\bar\eta)$.
If $m''_i\neq 0$, $l''_i(h^{1,2,3})=0$
because $l''_i$ is a combination of elements
of the form $l^{\bar\lambda(x),\bar\mu(y)}_{\bf c,\bf d}$ and $l^{\bar\nu(xy)}_{\bf e}$
for some pairs $(\bf c,\bar\lambda)$, $(\bf d,\bar\mu)$ and $(\bf e,\bar\nu)$.
Since $\delta(l^{\bar\zeta(x),\bar\eta(y)}_{\bf a,\bf b})(1\otimes h^{1,23,4}h^{1,2,3})=
\delta(l^{\bar\zeta(x),\bar\eta(y)}_{\bf a,\bf b})(h^{1,23,4}\otimes h^{1,2,3})$,
$$l^{\bar\zeta(x),\bar\eta(y)}_{\bf a,\bf b}(h^{1,23,4}h^{1,2,3})=
\sum_i l'_i(h^{1,23,4})\otimes l''_i(h^{1,2,3})
=l^{\bar\zeta(x),\bar\eta(y)}_{\bf a,\bf b}(h^{1,23,4})
=l^{\bar\zeta\bar\eta}_{\bf a\bf b}(h).$$
For the second equality we use the assumption $c_0(h)=1$.
\qed
\end{pf}

\begin{lem}\label{lemma4}
Let $(g,h)\in U \frak F_2\times U \frak F_{N+1}$ 
be a pair satisfying $c_0(h)=1$, $c_{A}(h)=0$ and \eqref{mixed pentagon}.
Suppose that $(\bf a,\bar\zeta)$ and $(\bf b,\bar\eta)$ are admissible. Then
$$l^{\bar\eta(y),\bar\zeta(x)}_{\bf b,\bf a}(h^{1,23,4}h^{1,2,3})=
l^{\bar\eta\bar\zeta}_{\bf b\bf a}(h).$$
\end{lem}

\begin{pf}
It follows $c_0(g)=1$ by our assumptions $c_0(h)=1$ and \eqref{mixed pentagon}.
Consider the embedding of Hopf algebra
$i_{2,3,4}:U\frak F_{2}\hookrightarrow U\frak t^0_{4,N}$
sending $A\mapsto t^{23}(0)$ and $B\mapsto t^{34}(0)$
(geometrically caused by the exceptional divisor obtained by blowing up at $(x,y)=(1,1)$.)
Since $(i_{2,3,4}\otimes id)(\Exp\Omega^{(N)}_4)
=\Exp\Omega^{(N)}_4(z)^{2,3,4}
\in U\frak t^0_{4,N}\widehat\otimes_{\bf Q} H^0{\bar B}({\mathcal M}_{0,4})
\simeq H^0{\bar B}({\mathcal M}^{(N)}_{0,5})^*\widehat\otimes_{\bf C} H^0{\bar B}({\mathcal M}_{0,4})$,
it induces the morphism
$$i^*_{2,3,4}:H^0{\bar B}({\mathcal M}^{(N)}_{0,5})\to H^0{\bar B}({\mathcal M}_{0,4})$$ 
which gives 
$i^*_{2,3,4}([\frac{dx}{x}])=0$,
$i^*_{2,3,4}([\frac{dx}{1-x}])=[\frac{dz}{z}]$,
$i^*_{2,3,4}([\frac{dx}{\zeta_N^a-x}])=0$ ($a\neq 0$),
$i^*_{2,3,4}([\frac{dy}{y}])=0$,
$i^*_{2,3,4}([\frac{dy}{1-y}])=[\frac{dz}{1-z}]$ and
$i^*_{2,3,4}([\frac{dy}{\zeta_N^a-y}])=0$ ($a\neq 0$),
$i^*_{2,3,4}([\frac{xdy+ydx}{\zeta^a_N-xy}])=0$ ($a\in{\bf Z}/N{\bf Z}$).
In each term of the expression 
$l^{\bar\eta(y),\bar\zeta(x)}_{\bf b,\bf a}
=\sum_{I=(i_m,\cdots,i_1)}c_I[\omega_{i_m}|\cdots|\omega_{i_1}]$,
the first component $\omega_{i_m}$ is always one of 
$\frac{dx}{x}$, $\frac{dy}{y}$, $\frac{dx}{\zeta_N^a-x}$ 
and $\frac{dy}{\zeta_N^a-y}$ for $a\neq 0$
because both $({\bf a},\bar\zeta)$ and $({\bf b},\bar\eta)$ are admissible.
So $i_{2,3,4}^*(l'_i)=0$ unless $m'_i=0$.
Therefore
$$l^{\bar\eta(y),\bar\zeta(x)}_{\bf b,\bf a}(g^{2,3,4}h^{1,23,4}h^{1,2,3})=
\sum_i l'_i(g^{2,3,4})\otimes l''_i(h^{1,23,4}h^{1,2,3})
=l^{\bar\eta(y),\bar\zeta(x)}_{\bf b,\bf a}(h^{1,23,4}h^{1,2,3})$$
by $c_0(g)=1$.
So by our assumption,
$$l^{\bar\eta(y),\bar\zeta(x)}_{\bf b,\bf a}(h^{1,23,4}h^{1,2,3})
=l^{\bar\eta(y),\bar\zeta(x)}_{\bf b,\bf a}(g^{2,3,4}h^{1,23,4}h^{1,2,3})
=l^{\bar\eta(y),\bar\zeta(x)}_{\bf b,\bf a}(h^{1,2,34}h^{12,3,4}).$$

By the same arguments to the last two paragraphs of the proof of lemma \ref{lemma3},
\begin{equation}\label{12,3,4}
i^*_{12,3,4}(l^{\bar\eta(y),\bar\zeta(x)}_{\bf b,\bf a})=0, \qquad
i^*_{12,3,4}(l^{\bar\zeta(xy)}_{\bf a})=0,
\end{equation}
\begin{equation*}\label{1,2,34}
i^*_{1,2,34}(l^{\bar\eta(y),\bar\zeta(x)}_{\bf b,\bf a})
=l^{\bar\eta\bar\zeta}_{\bf b,\bf a}, \qquad
i^*_{1,2,34}(l^{\bar\zeta(xy)}_{\bf a})=l^{\bar\zeta}_{\bf a}
\end{equation*}
for admissible pairs
$({\bf a},\bar\zeta)$ and $({\bf b},\bar\eta)$,
from which we can deduce
$$l^{\bar\eta(y),\bar\zeta(x)}_{\bf b,\bf a}(h^{1,2,34}h^{12,3,4})
=l^{\bar\eta\bar\zeta}_{\bf b\bf a}(h).$$
\qed
\end{pf}

\begin{lem}\label{lemma5}
Let $h\in U \frak F_{N+1}$ with $c_0(h)=1$ and $c_{A}(h)=0$.
Then
$$l^{\bar\zeta(x)}_{\bf a}(e^{T\{t^{23}(0)+t^{24}(0)+t^{34}(0)\}}h^{1,23,4}h^{1,2,3})=l^{\bar\zeta,I}_{\bf a}(h),$$
$$l^{\bar\zeta(y)}_{\bf a}(e^{T\{t^{23}(0)+t^{24}(0)+t^{34}(0)\}}h^{1,23,4}h^{1,2,3})=l^{\bar\zeta,I}_{\bf a}(h),$$
$$l^{\bar\zeta(xy)}_{\bf a}(e^{T\{t^{23}(0)+t^{24}(0)+t^{34}(0)\}}h^{1,23,4}h^{1,2,3})=l^{\bar\zeta,I}_{\bf a}(h),$$
$$l^{\bar\zeta(x),\bar\eta(y)}_{\bf a,\bf b}(e^{T\{t^{23}(0)+t^{24}(0)+t^{34}(0)\}}h^{1,23,4}h^{1,2,3})
=l^{\bar\zeta \bar\eta,I}_{\bf a\bf b}(h)$$
for any pairs $(\bf a,\bar\zeta)$ and $(\bf b,\bar\eta)$.
\end{lem}

\begin{pf}
By the arguments in lemma \ref{lemma3} and our assumption $c_{A}(h)=0$,
\begin{align*}
l^{\bar\zeta(x)}_{\bf a}
(&e^{T\{t^{23}(0)+t^{24}(0)+t^{34}(0)\}}h^{1,23,4}h^{1,2,3})
=l^{\bar\zeta}_{\bf a}
(p_4(e^{T\{t^{23}(0)+t^{24}(0)+t^{34}(0)\}}h^{1,23,4}h^{1,2,3}))\\
=&l^{\bar\zeta}_{\bf a}(e^{TB(0)}h)=l^{\bar\zeta,I}_{\bf a}(h),
\end{align*}
\begin{align*}
l^{\bar\eta(y)}_{\bf a}
(&e^{T\{t^{23}(0)+t^{24}(0)+t^{34}(0)\}}h^{1,23,4}h^{1,2,3})
=l^{\bar\eta}_{\bf a}
(p_2(e^{T\{t^{23}(0)+t^{24}(0)+t^{34}(0)\}}h^{1,23,4}h^{1,2,3}))\\
=&l^{\bar\eta}_{\bf a}(e^{TB(0)}h)=l^{\bar\eta,I}_{\bf a}(h),
\end{align*}
\begin{align*}
l^{\bar\zeta(xy)}_{\bf a}
(&e^{T\{t^{23}(0)+t^{24}(0)+t^{34}(0)\}}h^{1,23,4}h^{1,2,3})
=l^{\bar\zeta}_{\bf a}
(p_3(e^{T\{t^{23}(0)+t^{24}(0)+t^{34}(0)\}}h^{1,23,4}h^{1,2,3}))\\
=&l^{\bar\zeta}_{\bf a}(e^{TB(0)}h)=l^{\bar\zeta,I}_{\bf a}(h).
\end{align*}
By $c_0(h)=1$,
\begin{align*}
l^{\bar\zeta(x),\bar\eta(y)}_{\bf a,\bf b}
&(e^{T\{t^{23}(0)+t^{24}(0)+t^{34}(0)\}}h^{1,23,4}h^{1,2,3})
=l^{\bar\zeta(x),\bar\eta(y)}_{\bf a,\bf b}
(e^{T\{t^{24}(0)+t^{34}(0)\}}h^{1,23,4})\\
=&l^{\bar\zeta \bar\eta}_{\bf a\bf b}(e^{TB(0)}h)
=l^{\bar\zeta \bar\eta,I}_{\bf a\bf b}(h).
\end{align*}
As for the last equation, we use the following trick:
\begin{align*}
e^{Tt^{23}(0)}&e^{T\{t^{24}(0)+t^{34}(0)\}}h^{1,23,4}h^{1,2,3}
=e^{T\{t^{24}(0)+t^{34}(0)\}}e^{Tt^{23}(0)}h^{1,23,4}h^{1,2,3} \\
&=e^{T\{t^{24}(0)+t^{34}(0)\}}h^{1,23,4}e^{Tt^{23}(0)}h^{1,2,3}.
\end{align*}
\qed
\end{pf}

\begin{lem}\label{lemma6}
Let $(g,h)\in U \frak F_2\times U \frak F_{N+1}$ 
be a pair satisfying $c_0(h)=1$, $c_A(h)=c_{B(0)^n}(h)=0$ for all $n\in\bf N$
and \eqref{mixed pentagon}.
Suppose that $(\bf a,\bar\zeta)$ is admissible. Then
$$l^{\bar\eta(y),\bar\zeta(x)}_{\bf b,\bf a}
(e^{T\{t^{23}(0)+t^{24}(0)+t^{34}(0)\}}h^{1,23,4}h^{1,2,3})
=l^{\bar\eta\bar\zeta}_{\bf b\bf a}(h).$$
\end{lem}

\begin{pf}
Express 
$\delta(l^{\bar\eta(y),\bar\zeta(x)}_{\bf b,\bf a})=\sum_i l'_i\otimes l''_i$
with $\deg l'_i={m'_i}$ and $\deg l''_i={m''_i}$ for some $m'_i$ and $m''_i$
such that $m'_i+m''_i=wt({\bf a},\bar\zeta)+wt({\bf b},\bar\eta)$.
Since  $(\bf a,\bar\zeta)$ is admissible,
$i^*_{2,3,4}(l'_i)$ is of the form $\alpha[\frac{dz}{1-z}|\cdots|\frac{dz}{1-z}]$
with $\alpha\in {\bf Q}$.
But by our assumption $c_{B(0)^n}(h)=0$,
$i^*_{2,3,4}(l'_i)=0$ unless $m'_i=0$.
Thus
\begin{align*}
l^{\bar\eta(y),\bar\zeta(x)}_{\bf b,\bf a}&
(e^{T\{t^{23}(0)+t^{24}(0)+t^{34}(0)\}}h^{1,23,4}h^{1,2,3}) \\
&=l^{\bar\eta(y),\bar\zeta(x)}_{\bf b,\bf a}
(g^{2,3,4}e^{T\{t^{23}(0)+t^{24}(0)+t^{34}(0)\}}h^{1,23,4}h^{1,2,3}).
\end{align*}
By \eqref{mixed pentagon},
\begin{align*}
g^{2,3,4} & e^{T\{t^{23}(0)+t^{24}(0)+t^{34}(0)\}}h^{1,23,4}h^{1,2,3} \\
&=e^{T\{t^{23}(0)+t^{24}(0)+t^{34}(0)\}}g^{2,3,4}h^{1,23,4}h^{1,2,3} \\
&=e^{T\{t^{23}(0)+t^{24}(0)+t^{34}(0)\}}h^{1,2,34}h^{12,3,4} \\
&=e^{T\{t^{23}(0)+t^{24}(0)\}}e^{Tt^{34}(0)}h^{1,2,34}h^{12,3,4} \\
&=e^{T\{t^{23}(0)+t^{24}(0)\}}h^{1,2,34}e^{Tt^{34}(0)}h^{12,3,4}.
\end{align*}
By \eqref{12,3,4},
\begin{align*}
l^{\bar\eta(y),\bar\zeta(x)}_{\bf b,\bf a}&
(e^{T\{t^{23}(0)+t^{24}(0)\}}h^{1,2,34}e^{Tt^{34}(0)}h^{12,3,4})
=l^{\bar\eta(y),\bar\zeta(x)}_{\bf b,\bf a}
(e^{T\{t^{23}(0)+t^{24}(0)\}}h^{1,2,34}) \\
&=l^{\bar\eta\bar\zeta}_{\bf b\bf a}(e^{TB(0)}h)
=l^{\bar\eta\bar\zeta,I}_{\bf b\bf a}(h)
=l^{\bar\eta\bar\zeta}_{\bf b\bf a}(h).
\end{align*}
The last equality follows from the admissibility of $(\bf a,\bar\zeta)$
\qed
\end{pf}



\begin{thebibliography}{Utah}

\bibitem[AK]{AK}
Arakawa, T., Kaneko, M., 
{\it On multiple $L$-values}, 
J. Math. Soc. Japan 56 (2004), no.4, 967--991. 

\bibitem[B]{B}
Broadhurst, D. J.,
{\it Massive 3-loop Feynman diagrams reducible to $\rm SC^*$ primitives of algebras of the sixth root of unity},
Eur.Phys.J.C Part.Fields 8(1999), no. 2, 313--333.

\bibitem[BK]{BK}
Broadhurst, D.J., Kreimer, D.,
{\it Association of multiple zeta values with positive knots
via Feynman diagrams up to 9 loops},
Physics Letters B, Vol 393, Issues 3-4, 1997, 403-412.

\bibitem[Br]{Br}
Brown, F.,
{\it Mixed Tate motives over $\bold Z$},
Annals of Math., volume 175, no. 1 (2012), 949-976.

\bibitem[C]{C}
Chen, K. T.,
{\it Iterated path integrals},  Bull. Amer. Math. Soc.  83  (1977), no. 5, 831--879.

\bibitem[De89]{De}
Deligne, P.,
{\it Le groupe fondamental de la droite projective moins trois points},
Galois groups over ${Q}$ (Berkeley, CA, 1987), 79--297, Math. S. Res. Inst. 
Publ., 16, Springer, New York-Berlin, 1989.

\bibitem[De08]{De08}
\bysame,
{\it Le groupe fondemental de ${\bf G}_m\backslash\mu_N$, pour $N=2,3,4,6$, ou $8$},
preprint, 2008, available from
www.math.ias.edu/files/deligne/121108Fondamental.pdf

\bibitem[DG]{DG}
\bysame, Goncharov, A. B.,
{\it Groupes fondamentaux motiviques de Tate mixte},
Ann. Sci. Ecole Norm. Sup. (4)  38  (2005),  no. 1, 1--56.

\bibitem[Dr]{Dr}
Drinfel'd, V. G.,
{\it On quasitriangular quasi-Hopf algebras and  
a group closely connected with ${\rm Gal}(\overline{Q}/{Q})$},
Leningrad Math. J. 2 (1991), no. 4, 829--860.

\bibitem[E]{E}
Enriquez, B.,
{\it Quasi-reflection algebras and cyclotomic associators},
Selecta Math. (N.S.) 13 (2007), no. 3, 391-463.

\bibitem[EF]{EF}
\bysame, Furusho, H.,
{\it Mixed Pentagon, octagon and Broadhurst duality equation},
Journal of Pure and Applied Algebra, Vol 216, Issue 4, (2012), 982-995. 

\bibitem[F1]{F03}
Furusho, H.,
{\it The multiple zeta value algebra and the stable derivation algebra}, 
Publ. Res. Inst. Math. Sci. Vol 39. no 4. (2003). 695-720.

\bibitem[F2]{F10}
\bysame,
{\it Pentagon and hexagon equations},
Annals of Mathematics, Vol. 171 (2010), No. 1, 545-556.

\bibitem[F3]{F08}
\bysame,
{\it Double shuffle relation for associators},
Annals of Mathematics, Vol. 174 (2011), No. 1, 341-360.

\bibitem[G]{G}
Goncharov, A. B.,
{\it The dihedral Lie algebras and Galois symmetries of
$\pi_1^{(l)}(\Bbb P^1-(\{0,\infty\}\cup\mu_N))$},
Duke Math. J. 110 (2001),  no. 3, 397--487.

\bibitem[K]{K}
Kohno, T.,
{\it Bar complex of the Orlik-Solomon algebra},
Arrangements in Boston: a Conference on Hyperplane Arrangements (1999),
Topology Appl. 118 (2002), no. 1-2, 147-157.


\bibitem[LM]{LM}
Le, T.T.Q., Murakami, J.;
{\it Kontsevich's integral for the Kauffman polynomial},
Nagoya Math. J. 142 (1996), 39--65.

\bibitem[R]{R}
Racinet, G.,
{\it Doubles melanges des polylogarithmes multiples aux racines de l'unite},
Publ. Math. Inst. Hautes Etudes Sci. No. 95 (2002), 185--231.

\bibitem[Z]{Z}
Zhao, Jianqiang.,
{\it Multiple polylogarithm values at roots of unity},
C. R. Math. Acad. Sci. Paris  346  (2008),  no. 19-20, 1029--1032. 

\end{thebibliography}
\end{document}